\documentclass[final,5p,times,twocolumn]{elsarticle}

\usepackage{amssymb}
\usepackage{amsmath} 
\usepackage{color}
\usepackage{algorithm}
\usepackage{algpseudocode}
\usepackage{float}
\usepackage{multirow}
\usepackage{caption}

\algnewcommand{\algorithmicforeach}{\textbccf{for each}}
\algdef{SE}[FOR]{ForEach}{EndForEach}[1]
  {\algorithmicforeach\ #1\ \algorithmicdo}
  {\algorithmicend\ \algorithmicforeach}

\usepackage[switch]{lineno}

\newcommand{\mb}{\mathbf}

\renewcommand{\arraystretch}{2} 

\begin{document}

\begin{frontmatter}



\title{Distortion-minimized de-homogenization for optimization of cell-size distribution in TPMS structures}


\author{Hiroki Kawabe\corref{cor1}}
\cortext[cor1]{Corresponding author.}
\ead{kawabe@syd.mech.eng.osaka-u.ac.jp}
\author{Kaito Ohtani}
\author{Yusibo Yang}
\author{Musaddiq Al Ali}
\author{Kentaro Yaji}


\address{Department of Mechanical Engineering, Graduate School of Engineering, The University of Osaka, 2-1, Yamadaoka, Suita, Osaka, 565-0871, Japan}

\begin{abstract}
This paper presents a homogenized topology optimization (TO) method for spatially optimizing cell-size distribution of triply-periodic minimal surface (TPMS) structures, with high accuracy in the optimized structural response after de-homogenization.
To achieve this, we introduce a novel de-homogenization technique that directly minimizes the difference between the wavenumbers obtained from the target and actual size distributions.
This minimization problem is efficiently solved as a typical Poisson's equation utilizing the discrete cosine transform (DCT).
We first verify the proposed de-homogenization method through numerical examples, showcasing its capability in significantly reducing the known distortion of the de-homogenized TPMS structures from the conventional periodic modulation (PM) method.
Then, we apply the proposed method to a stiffness maximization problem, to demonstrate its effectiveness in improving the structural response compared to the PM method, through numerical examples and experimental validation.
The proposed method successfully reduced the distortion of the de-homogenized structures compared to the PM method, leading to 0.8\% difference in the strain energy compared to the homogenized model, as opposed to 63.6\% difference in the PM method.
The optimized structure from the proposed method shows a significant improvement in the strain energy by 50.1\% compared to the uniform case in the FE analysis on the de-homogenized models, while the PM method results in a significant decrease of 45.8\%.
The experimental validation shows that the effective stiffness of the optimized structure from the proposed method is 54.2\% higher than that of the uniform case, while the PM method results in a significant decrease by 77.3\%.
These results exhibit the proposed method effectively increases the accuracy of the de-homogenization, thereby maximizing the potential of the homogenized TO for the spatial cell-size optimization of TPMS structures.
\end{abstract}

%


\end{frontmatter}



\section{Introduction}
\label{sec:intr}

Porous structures~\cite{benedetti2021architected}, such as truss lattice structures~\cite{pan2020design} and honeycomb structures~\cite{qi2021advanced}, are well known for their lightweight nature combined with high stiffness, high strength, and high surface area.
Owing to these advantageous properties, these structures have been widely used in various fields, including structural components~\cite{khan2024systematic}, medical applications~\cite{yuan2019additive}, and thermal systems~\cite{sajjad2022manufacturing}.
These traditional porous structures can be easily manufactured using conventional manufacturing methods, such as casting and machining.
However, they inherently possess sharp corners and edges, which can lead to stress concentration and potential failure under mechanical loading~\cite{zhao2022design}.

With the advancement of additive manufacturing (AM), the fabrication of complex structures has become increasingly feasible, enabling the development of highly efficient porous architectures with enhanced performance~\cite{mohamed2022enabling,shimoda2023micropore}.
Recently, triply-periodic minimal surface (TPMS) structures~\cite{alketan2019multifunctional} have gained much attention as a promising alternative to traditional porous structures, due to their smooth surfaces without sharp corners and high surface area to volume ratio, leading to high specific stiffness and strength~\cite{feng2022triply}.
TPMS structures are represented by implicit functions, which allow for easy control of the thickness and size of the TPMS structure, full connectivity, and quasi-self-supporting nature.


Utilizing these capabilities of TPMS structures, various design methods have been proposed to generate functionally graded TPMS structures with spatially varying thickness distributions~\cite{qiu2023mechanical,jia2020experimental} and size distributions~\cite{yan2020strong, hu2019lightweight}, or combining multiple types of TPMS structures~\cite{yoo2015advanced, maskery2018effective}.
Due to its simplicity, the thickness variation method has been widely used in various applications, such as the design of bone scaffolds~\cite{cohen2025optimal}, tissue scaffolds~\cite{melchels2010mathematically}, heat exchangers~\cite{chen2023heat}, and structural components~\cite{yu2019investigation}.
Despite its effectiveness in enhancing mechanical performance of the TPMS structures, the enhancement is limited because their shape is not drastically changed.
The multiple type method can offer more drastic shape changes, by mixing multiple types of TPMS structures~\cite{alketan2020functionally,wang2022crashworthiness,xi2023multi,ning2025high}.
However, the multiple type method requires a transition region with special treatment like thickening, which may compromise the performance of the structure~\cite{ren2021transition}.
In contrast, the size variation method can offer more gradual changes in the TPMS structure without the need for a distinct transition region, potentially preserving the mechanical performance of the structure.
The method has been used in various applications, such as the design of medical implants~\cite{liu2020additively}, heat exchangers~\cite{gao2023influence}, and structural components~\cite{alketan2020functionally, caiazzo2022metal}.
However, the design of the size distribution has been mostly based on heuristic approaches, with the number of design variables often limited to a few.

To overcome these limitations, topology optimization (TO)~\cite{bendsoe2003topology} can be promising.
TO creates innovative and high-performance designs based on the mathematical formulation of the optimization problem and numerical optimization algorithms, which enables the efficient exploration of a vast design space without the need for human intuition.
However, it is inherently challenging to utilize TO for designing complex geometries such as TPMS structures, due to the high-resolution discretization required to accurately capture their complex geometries, leading to extremely high computational cost.

To address this issue, Hu et al.~\cite{hu2022efficient} adopted a super-element technique~\cite{nguyen2010computational}, which divides the design domain into coarse subdomains called super-elements, while each super-element is divided into fine elements called background elements.
This method can reduce the computational cost by performing most of the finite element analysis (FEA) processes on the super-element level, while the design variables such as parameters for controlling the thickness and size of the TPMS structures are defined on the background element level.
The super-element-based TO proved to be effective in designing high-stiffness functionally graded TPMS structures, later extended to the heat dissipation problem~\cite{wang2022efficient}.
Ning et al.~\cite{ning2026data} extended the super-element-based TO to the design of hybrid TPMS structures by introducing a data-driven surrogate model.
These methods, however, require the integral computation of the stiffness matrix for each super-element at each iteration, which can be computationally prohibitive especially when each super-element needs high-resolution to capture small cell sizes or drastic changes in geometry.

Another approach to tackle the computational cost issue is to use a homogenization method~\cite{bendsoe1988generating}.
This homogenization-based TO method assumes that the effective properties of the TPMS structure can be estimated by analytical or numerical homogenization methods, utilizing its periodic nature.
By optimizing the distribution of the local structural features in the homogenized model and reconstructing the actual geometry of the TPMS structure from the optimized distribution of structural features, the method can efficiently design functionally graded TPMS structures with a large number of design variables.
This reconstruction process is called de-homogenization, which is a critical step to ensure that the optimized design from the homogenized model can be accurately reflected in the reconstructed TPMS structure to achieve the optimized performance.
Li et al.~\cite{li2019design} proposed an asymptotic homogenization-based TO method for the stiffness problem and the thermal conductivity problem, optimizing the relative density distribution of the TPMS structures.
Str{\"o}mberg~\cite{stromberg2021optimal} utilized a numerical homogenization method using three repeated cells for the stiffness problem, optimizing the relative density distribution.
Feng et al.~\cite{feng2022stiffness} proposed an energy-based numerical homogenization method with periodic boundary conditions for the stiffness problem, optimizing multiple TPMS cell types and relative density distribution.
Zhao et al.~\cite{zhao2023tpms} proposed a method for the stiffness problem, optimizing the relative density distribution of interpenetrating TPMS structures with multiple cell types mixed.
Str{\"o}mberg~\cite{stromberg2024new} extended their method~\cite{stromberg2021optimal} to additionally optimize void and lattice-infill regions.
The homogenization-based TO has also been successfully applied to thermal-fluid problems, such as the design of heat sinks~\cite{men2025topology}.


%

Despite the success of the homogenization-based TO in designing functionally graded TPMS structures, the method has a critical issue in the de-homogenization process when the size distribution is spatially varied.
In the de-homogenization process, the actual geometry of the TPMS structure is generated by evaluating the implicit function of the TPMS with the phase distributions obtained from the optimized size distribution.
When calculating the phase distributions, a traditional method, called periodic modulation (PM) method, directly uses the optimized size distribution to define the phase distributions.
The PM method inherently causes distortion of the de-homogenized TPMS structure due to the discrepancy between the gradients of the phase distributions, which represent the real wavenumbers corresponding to the actual size distribution, and the desired wavenumbers obtained from the optimized size distribution in the homogenization-based TO.
This undesired distortion often leads to a significant gap between the performance of the optimized design and the actual performance of the de-homogenized structure.

In the context of the functionally graded TPMS structures in general, a few studies have proposed methods to mitigate this issue, although none of them have been applied to the homogenization-based TO for size optimization of TPMS structures.
Liu et al.~\cite{liu2018functionally} introduced a constraint on the gradients of the phase distributions to keep the structural consistency of the TPMS structure, which can avoid severe distortion but limits the design freedom and the performance of the optimized structure.
Tian et al.~\cite{tian2024continuous} proposed a method to optimize the phase distributions by minimizing an indirect measure of the distortion of the TPMS structure, enforcing orthogonality among the gradients of the phase distributions.
This phase optimization is effective in ensuring the local orthogonality of each phase distribution and thus reducing the distortion to some extent, but still results in significant distortion when the TPMS size varies drastically in space after the homogenization-based TO.
Wang and Zhong~\cite{wang2025optimization} proposed a method to optimize the phase distributions by directly minimizing the discrepancy between the gradients of the phase distributions and the desired wavenumbers of the TPMS structure.
However, their method is based on Legendre polynomials to reduce the number of design variables and thus the computational cost, which still suffers from high computational cost especially under high-order polynomials.
Moreover, the method requires a three-step preprocessing of the target size distribution, including discretization, Gaussian smoothing, and remapping.
This not only necessitates careful tuning of the preprocessing parameters in addition to the Legendre polynomial order, but also alters the size distribution from the optimized one obtained by the homogenization-based TO.
Therefore, a more effective and efficient method to optimize the phase distributions is needed to minimize the distortion of the de-homogenized TPMS structure and thus to better reflect the performance of the optimized design in the de-homogenized structure.

To achieve this, we propose a method to directly minimize the discrepancy between the gradients of the phase distributions and the desired wavenumbers obtained from the size distribution optimized by the homogenization-based TO, by solving the least-squares problem separately for each phase distribution.
Each least-squares problem is formulated as a Poisson's equation with Neumann boundary conditions, which can be solved efficiently without the need for any preprocessing or approximation of the phase distributions obtained from the optimized size distribution.
We first evaluate the effectiveness of the proposed method in minimizing the distortion of the de-homogenized TPMS structure under various conditions of the size distribution, including a simple unidirectional size variation and more complex 3D cases, in comparison with the conventional PM method.
Then we apply the proposed method to the de-homogenization of the size distribution optimized by the homogenization-based TO for a stiffness maximization problem. 
This result demonstrates the effectiveness of the proposed method in closely reflecting the optimized size and thus significantly improving the performance of the de-homogenized TPMS structure compared to those of the conventional PM method and non-optimized size distribution.
Finally, we fabricate those de-homogenized TPMS structures using 3D printing and experimentally measure their stiffness under a loading condition that replicates the numerical simulation, to validate the effectiveness proved in the numerical examples.

The rest of this paper is organized as follows.
Section~\ref{sec:mthd} describes the methodology of the proposed method, including the formulation of the homogenization-based TO for size optimization and the proposed de-homogenization method.
Section~\ref{sec:impl} details the numerical implementation of the proposed method, including the numerical examples to evaluate the effectiveness of the proposed method.
Section~\ref{sec:num} presents the results of the homogenization-based TO and the de-homogenization with the proposed method, whereas Section~\ref{sec:exp} describes the experimental setup and results of the stiffness measurement of the 3D printed TPMS structures.
Finally, Section~\ref{sec:con} concludes the paper with a summary of the findings and future work.

\section{Methodology}
\label{sec:mthd}
\begin{figure}[t]
  \centering
  \includegraphics[width=\columnwidth]{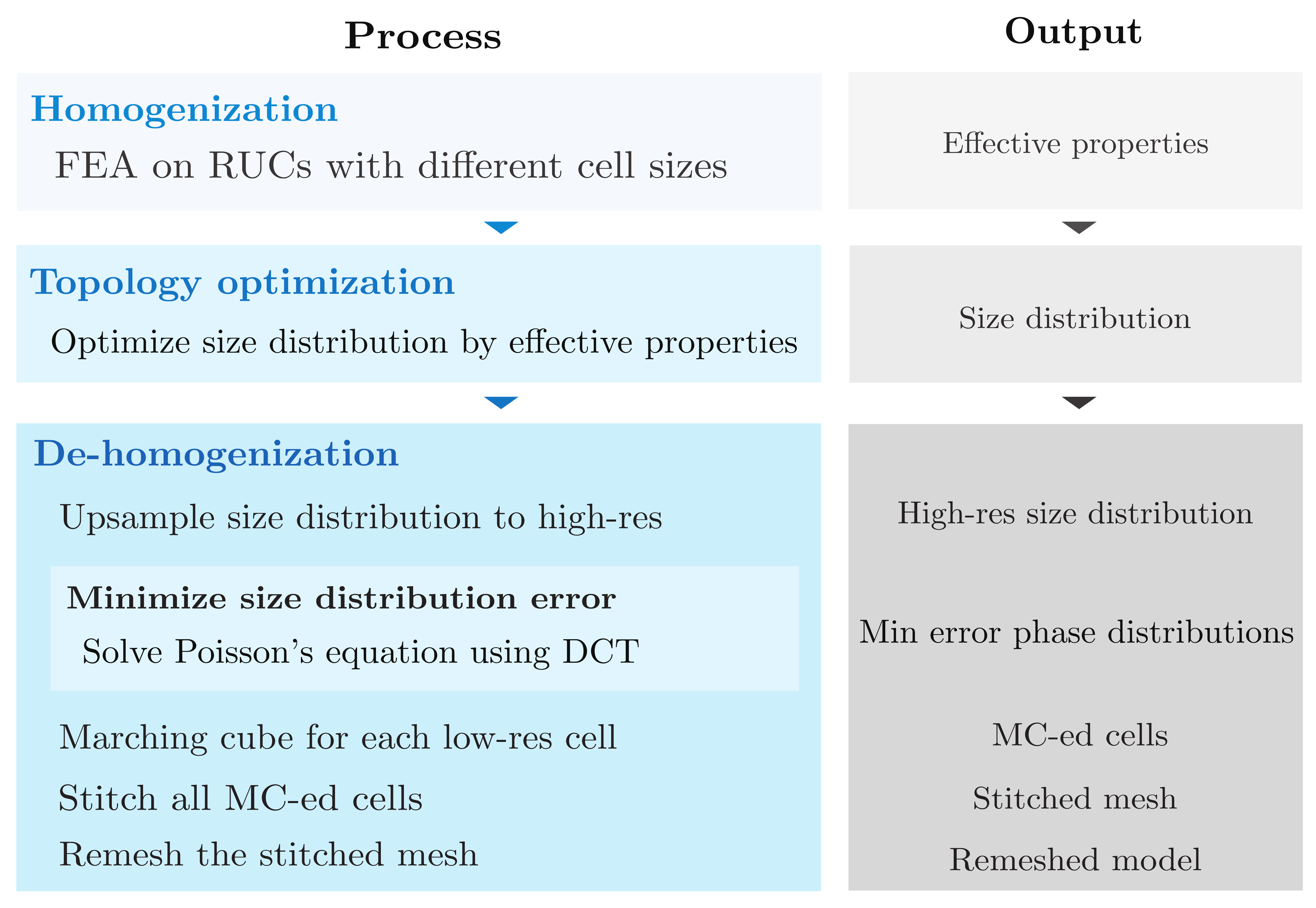}
  \caption{Overall framework of the proposed method.}
  \label{fig:framework}
\end{figure}

In order to efficiently design functionally graded TPMS structures with spatially varying size distributions, we adopt a homogenization-based TO method.
This method consists of the following steps: 1) effective property estimation by homogenization, 2) topology optimization of the size distribution using the estimated effective properties, and 3) de-homogenization with the optimized size distribution.
To minimize the distortion of the de-homogenized TPMS structure and thus to de-homogenize the optimized size distribution as accurately as possible, we introduce a method to optimize the phase distributions in the de-homogenization process after the homogenization-based TO.
Fig.~\ref{fig:framework} illustrates the overall framework of the proposed method.

\subsection{Size graded TPMS structures}
In general, TPMS is defined as the implicit surface that satisfies the following equation:
\begin{equation}
  F(\phi_x(\mb{r}), \phi_y(\mb{r}), \phi_z(\mb{r})) = 0
  \label{eq:tpms_general}
\end{equation}
\noindent where $F$ is the implicit function that defines the TPMS, and $\phi_x$, $\phi_y$, and $\phi_z$ are the phase distributions in the $x$-, $y$-, and $z$-directions, respectively; $\mb{r}=(x,y,z)$ is the position vector in the Cartesian coordinate system.

For example, the implicit function of the Schwarz P TPMS and Gyroid TPMS are given as follows:
\begin{equation}
  \begin{aligned}
    F_P (\mb{r}) & = \cos (\phi_x(\mb{r})) + \cos (\phi_y(\mb{r})) + \cos (\phi_z(\mb{r})) = 0, \\  
    F_G (\mb{r}) & = \sin (\phi_x(\mb{r})) \cos (\phi_y(\mb{r})) + \sin (\phi_y(\mb{r})) \cos (\phi_z(\mb{r})) \\
    & \ + \sin (\phi_z(\mb{r})) \cos (\phi_x(\mb{r})) = 0.
  \end{aligned}
  \label{eq:type}
\end{equation}
\noindent where $F_P$ and $F_G$ are the implicit functions of the Schwarz P TPMS and Gyroid TPMS, respectively. 
Then solid regions of the TPMS sheet with a certain thickness can be defined for each type as follows:
\begin{equation}
  \Omega_\text{solid} = \left\{\mb{r} \mid -\frac{c}{P} \leq F_\text{type} \leq \frac{c}{P} \right\},
  \label{eq:tpms_solid}
\end{equation}
\noindent where $c$ is the level parameter that controls the thickness of the TPMS sheet, and $P$ is the period that equals the size of the TPMS structure.

In the traditional PM method for controlling the size of the TPMS structure, the phase distributions are defined based on the spatially varying size $P(\mb{r})$ as follows:
\begin{equation}
  \phi_x^\text{PM}(\mb{r}) = \frac{2\pi}{P(\mb{r})} x, \quad
  \phi_y^\text{PM}(\mb{r}) = \frac{2\pi}{P(\mb{r})} y, \quad
  \phi_z^\text{PM}(\mb{r}) = \frac{2\pi}{P(\mb{r})} z,
  \label{eq:phase_traditional}
\end{equation}
\noindent where $\phi_s^\text{PM}$ ($s=x,y,z$) are the phase distributions from the PM method.
The spatially varying size $P(\mb{r})$ is originally intended to directly represent the local size of the TPMS structure.

\begin{figure}[t]
  \centering
  \includegraphics[width=\columnwidth]{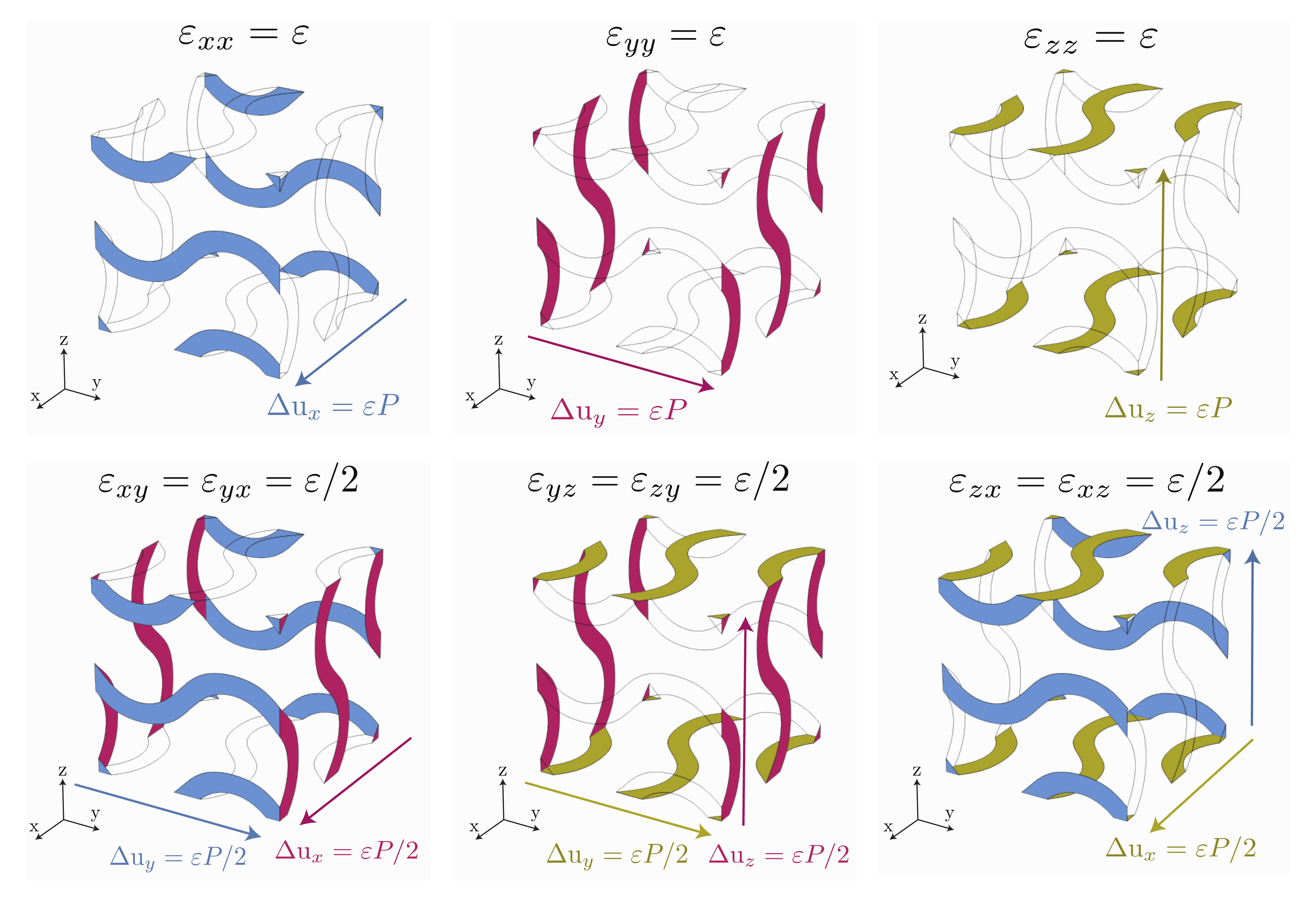}
  \caption{RUC with periodic boundary conditions applied to the opposite faces.}
  \label{fig:unit_cell}
  \end{figure}

\subsection{Distortion-minimized de-homogenization}
\label{subsec:dehom}

In Eq.~\ref{eq:phase_traditional}, when the size $P(\mb{r})$ is uniform over the entire design domain, the generated TPMS structure has a uniform size of $P$ and a constant wavenumber of $\omega = 2\pi / P$ in each direction.
This wavenumber can be calculated as the gradient of the phase distributions as follows:
\begin{equation}
  \nabla \phi_s^\text{Uni} = \omega \mb{e}_s, \quad
  \phi_s^\text{Uni} = \omega s, \quad
  s=x,y,z,
  \label{eq:grad_phase_traditional}
\end{equation}
\noindent where $\nabla$ is the gradient operator, $\mb{e}_s$ is the unit vector in the $s$-direction, and $\phi_s^\text{Uni}$ are the phase distributions for the uniform size case.

However, when the size $P(\mb{r})$ varies spatially, the real wavenumbers from the PM method are calculated as
\begin{equation}
  \nabla \phi_s^\text{PM} = \omega \mb{e}_s + s \nabla \omega,
  \label{eq:grad_phase_spatial}
\end{equation}
\noindent which differ from the intended wavenumber $\omega \mb{e}_s$ because of the additional terms including $\nabla \omega$.
This discrepancy causes distortion of the generated TPMS structure especially when the spatial variation of $P(\mb{r})$ is large, and/or when the distance from the origin is large.
To overcome this issue, we propose a method to optimize the phase distributions $\phi_s$ so that their gradients match the desired wavenumbers $\omega \mb{e}_s$ as closely as possible.
Note that the size $P(\mb{r})$ is optimized in the homogenization-based TO, and this method aims to find the phase distributions $\phi_s$ that best represent the optimized size distribution.
For this purpose, we formulate the following optimization problem for each phase distribution $\phi_s$:

\begin{equation}
  \min_{\phi_s} \quad \int_{\Omega} || \nabla \phi_s - \omega \mb{e}_s||^2 \, d\Omega
  \label{eq:opt_phase}
\end{equation}
\noindent where $\Omega$ is the design domain.
The Euler-Lagrange equation of the problem~\eqref{eq:opt_phase} is given as
\begin{equation}
  \begin{aligned}
    \Delta \phi_s & = \nabla \cdot \boldsymbol{\omega}_s \quad \text{in} \quad \Omega, \\
    \frac{\partial \phi_s}{\partial n} & = \boldsymbol{\omega}_s \cdot \boldsymbol{n} \quad \text{on} \quad \partial \Omega,
  \end{aligned}
  \label{eq:euler_lagrange}
\end{equation}
\noindent where $\Delta$ is the Laplacian operator, $\boldsymbol{n}$ is the outward unit normal vector on the boundary $\partial \Omega$ of the design domain, and $\partial \phi_s / \partial n$ is the normal derivative of $\phi_s$ on $\partial \Omega$.
By solving Eq.~\eqref{eq:euler_lagrange} for each $s=x,y,z$, we can obtain the phase distributions $\phi_s$ that minimize the discrepancy between their gradients and the desired wavenumbers.
These optimized phase distributions can then be used in Eq.~\eqref{eq:tpms_general} to generate the TPMS structure that accurately reflects the spatially varying size $P(\mb{r})$.


\begin{figure}[t]
  \centering
  \includegraphics[width=0.9\columnwidth]{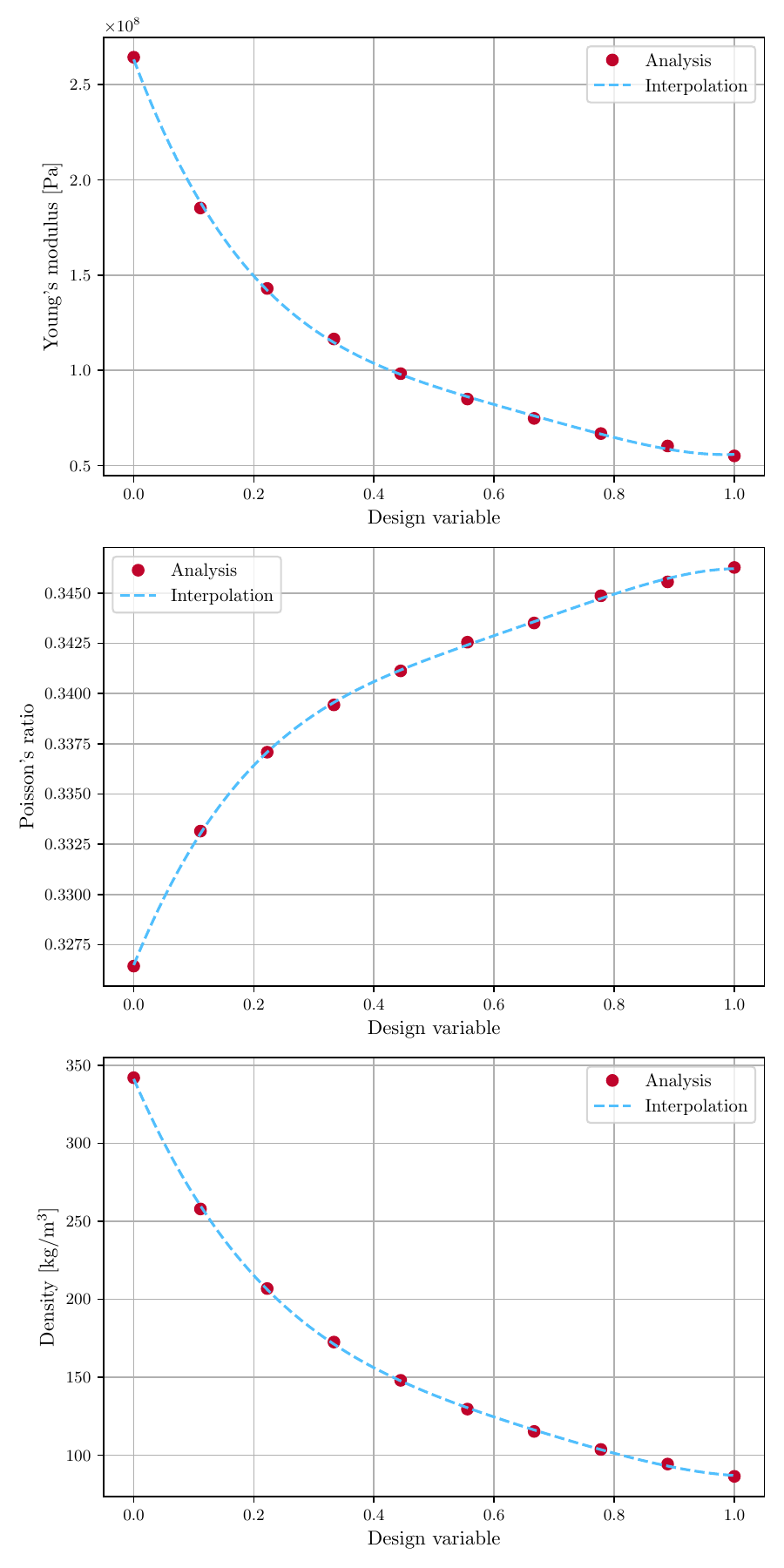}
  \caption{Interpolation of the effective properties of the TPMS structure based on the RUC analysis: Young's modulus (top), Poisson's ratio (center), and density (bottom) as functions of the design variable.}
  \label{fig:interpolation}
  \end{figure}

\begin{figure*}[!tb]
  \centering
  \includegraphics[width=0.8\textwidth]{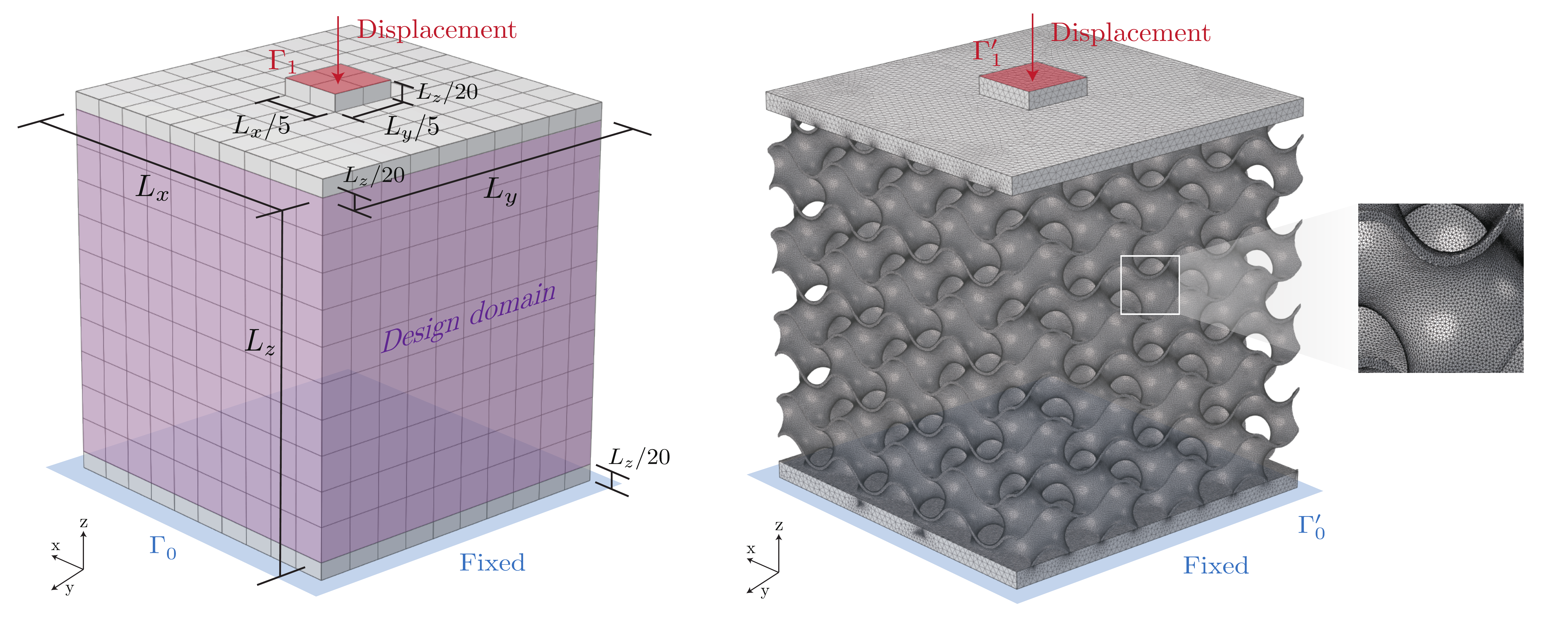}
  \caption{Simulation setup: the analysis domain and boundary conditions for the FEA of the homogenized model (left) and the de-homogenized model (right).}
  \label{fig:simulation_setup}
  \end{figure*}

\subsection{Effective property estimation by homogenization}
\label{subsec:hom}

In the homogenization process, we assume that the TPMS structure is locally periodic even when the size distribution is spatially varied during the optimization process.
Based on this assumption, the effective properties of the TPMS structure can be estimated by numerical homogenization methods, using repeating unit cells (RUCs) with periodic boundary conditions~\cite{drago2007micromacromechanical}.
The RUC is defined as a cube with a constant size of $P^{(k)} (k=1,2,\dots,N_\text{RUC})$ in each direction, following the phase distributions defined in Eq.~\eqref{eq:phase_traditional}.
The periodic boundary conditions are applied to each pair of opposite faces of the RUC, as follows:
\begin{equation}
  \mb{u}^{(k)}(\mb{r} + P^{(k)} \mb{e}_s) = \mb{u}^{(k)}(\mb{r}) + P^{(k)} \boldsymbol{\varepsilon}^{(k)} \mb{e}_s,
  \label{eq:periodic_bc}
\end{equation}
\noindent where $\mb{u}^{(k)}$ is the displacement vector for the $k$-th RUC, and $\boldsymbol{\varepsilon}^{(k)}$ is the macroscopic strain tensor for the $k$-th RUC.
Fig.~\ref{fig:unit_cell} illustrates the RUC with the periodic boundary conditions applied to the opposite faces for Gyroid TPMS.

The homogenized elasticity tensor $\mb{D}^{(k)}$ can be obtained from the following constitutive relation:
\begin{equation}
  \boldsymbol{\sigma}^{(k)} = \mb{D}^{(k)} : \boldsymbol{\varepsilon}^{(k)},
  \label{eq:stiffness_tensor}
\end{equation}
\noindent where $\boldsymbol{\sigma}^{(k)}$ is the macroscopic stress tensor for the $k$-th RUC.
FEA are performed on the RUC with the periodic boundary conditions under six independent loading conditions, each defined in Voigt notation: $\boldsymbol{\varepsilon}^{(k)}_1=[\varepsilon,0,0,0,0,0]^T$, $\boldsymbol{\varepsilon}^{(k)}_2=[0,\varepsilon,0,0,0,0]^T$, ..., $\boldsymbol{\varepsilon}^{(k)}_6=[0,0,0,0,0,\varepsilon]^T$ with a strain magnitude of $\varepsilon$.
For each loading condition, the macroscopic stress tensor $\boldsymbol{\sigma}^{(k)}$ is calculated as the volume average of the local stress field $\boldsymbol{\sigma}^\text{RUC}$ in the RUC, as follows:
\begin{equation}
  \boldsymbol{\sigma}^{(k)}_l = \frac{1}{{P^{(k)}}^3} \int_{\Omega_s} \boldsymbol{\sigma}^\text{RUC}(\mb{r}) \, d\Omega,
  \label{eq:macroscopic_stress}
\end{equation}
\noindent where $\Omega_s$ is the solid region of the RUC, and $l=1,2,\dots,6$ is the index of the loading condition.
By combining the calculated macroscopic stress tensors $\boldsymbol{\sigma}^{(k)}_l$ for the six loading conditions, we can solve for the components of the stiffness tensor $\mb{D}^{(k)}$, as follows:
\begin{equation}
  \mb{D}^{(k)} = \frac{1}{\varepsilon}\left[\boldsymbol{\sigma}^{(k)}_1 \ \boldsymbol{\sigma}^{(k)}_2 \ \boldsymbol{\sigma}^{(k)}_3 \ \boldsymbol{\sigma}^{(k)}_4 \ \boldsymbol{\sigma}^{(k)}_5 \ \boldsymbol{\sigma}^{(k)}_6 \right].
  \label{eq:stiffness_tensor_calculation}
\end{equation}

The compliance tensor $\mb{C}^{(k)}$ can be obtained by inverting the stiffness tensor $\mb{D}^{(k)}$.
This study assumes that the TPMS structure is isotropic, so the effective Young's modulus $E^{(k)}$ and the effective Poisson's ratio $\nu^{(k)}$ can be calculated from the two independent components of the compliance tensor $\mb{C}^{(k)}$ as follows:
\begin{equation}
  E^{(k)} = \frac{1}{S^{(k)}_{11}}, \quad
  \nu^{(k)} = -\frac{S^{(k)}_{12}}{S^{(k)}_{11}}.
  \label{eq:effective_properties}
\end{equation}
\noindent $S^{(k)}_{11}$ and $S^{(k)}_{12}$ are the components of the compliance tensor $\mb{C}^{(k)}$ in the Voigt notation.
The effective density $\rho^{(k)}$ can be calculated as the ratio of the volume of the solid region to the total volume of the RUC ${P^{(k)}}^3$, which is easily obtained from the FE model of the RUC.

By performing the homogenization process for each RUC with different sizes $P^{(k)}$, we can obtain the effective properties of the TPMS structure as interpolation functions of a normalized size $\bar{P}=(P-P_\text{min})/(P_\text{max}-P_\text{min})$, which can be used in the subsequent TO process, as follows:
\begin{equation}
  E^* = f_E(\bar{P}), \quad
  \nu^* = f_\nu(\bar{P}), \quad
  \rho^* = f_\rho(\bar{P}),
  \label{eq:effective_properties_functions}
\end{equation}
\noindent where $f_E$, $f_\nu$, and $f_\rho$ are the interpolation functions for the effective Young's modulus, effective Poisson's ratio, and effective density, respectively.

\subsection{Topology optimization}
\label{subsec:top}
In the homogenization-based TO, to optimize structural features in a continuous manner, we introduce a continuous design variable field $\gamma \in [0, 1]$ that represents the normalized value of the structural features within the design domain.
Here, the structural features refer to the size of the TPMS structure, which is expressed as a function of the design variable $\gamma$ as follows:
\begin{equation}
  P(\mb{r}) = P_\text{min} + \tilde{\gamma}(\mb{r}) (P_\text{max} - P_\text{min}),
  \label{eq:size_design_variable}
\end{equation}
\noindent where $P_\text{min}$ and $P_\text{max}$ are the minimum and maximum values of the size of the TPMS structure, respectively; $\tilde{\gamma}$ is the filtered design variable field obtained by applying a filter~\cite{bourdin2001filters} to the original design variable field $\gamma$ to mitigate overly intense local changes in the size distribution.
This filter averages the design variables within a certain neighborhood of each node, as follows:
\begin{align}
  w(\mathbf{r},\boldsymbol{\xi}) =
  \begin{cases}
    \dfrac{r_{\text{f}} - \lVert \mathbf{r}-\boldsymbol{\xi} \rVert}{r_{\text{f}}}, & \text{if } \lVert \mathbf{r}-\boldsymbol{\xi} \rVert < r_{\text{f}}, \\[6pt]
    \quad \quad \ 0, & \text{otherwise}.
  \end{cases}
\label{eq:weight}
\end{align}

\begin{align}
  \tilde{\gamma}(\mathbf{r})
  =
  \frac{\displaystyle
    \int_{\Omega} w(\mathbf{r},\boldsymbol{\xi})\,
                   \gamma(\boldsymbol{\xi})\,\mathrm{d}\boldsymbol{\xi}}
       {\displaystyle
    \int_{\Omega} w(\mathbf{r},\boldsymbol{\xi})\,
                   \mathrm{d}\boldsymbol{\xi}},
  \label{eq:filter_design_variable}
\end{align}
\noindent where $r_{\text{f}}$ is the filter radius parameter; $w(\mb{r}, \boldsymbol{\xi})$ is the weight function defined as a hat-type kernel.


The governing equations are formulated based on the effective properties obtained from the homogenization process, as follows:
\begin{equation}
\begin{aligned}
&  \nabla \cdot \boldsymbol{\sigma} = 0, \\
& \boldsymbol{\sigma} = \mb{D}^*(f_E(P), f_\nu(P)) : \boldsymbol{\varepsilon}, \\
& \boldsymbol{\varepsilon} = \frac{1}{2} \left( \nabla \mb{u} + (\nabla \mb{u})^\text{T} \right),
\label{eq:stiffness_equation}
\end{aligned}
\end{equation}
\noindent
where $\boldsymbol{\sigma}$ is the stress tensor, $\boldsymbol{\varepsilon}$ is the strain tensor, $\mb{u}$ is the displacement vector, and $\mb{D}^*$ is the effective elasticity tensor for isotropic materials, defined as follows:
{
\setlength{\arraycolsep}{1pt} 
\begin{equation}
  \begin{aligned}
    \mb{D}^* & = \frac{E^*}{(1+\nu^*)(1-2\nu^*)} \bar{\mb{D}}^*, \\
    \bar{\mb{D}}^* & = \begin{bmatrix}
  1-\nu^* & \nu^* & \nu^* & 0 & 0 & 0 \\
  \nu^* & 1-\nu^* & \nu^* & 0 & 0 & 0 \\
  \nu^* & \nu^* & 1-\nu^* & 0 & 0 & 0 \\
  0 & 0 & 0 & \frac{1-2\nu^*}{2} & 0 & 0 \\
  0 & 0 & 0 & 0 & \frac{1-2\nu^*}{2} & 0 \\
  0 & 0 & 0 & 0 & 0 & \frac{1-2\nu^*}{2}
  \end{bmatrix}.
  \end{aligned}
\label{eq:effective_elasticity_tensor}
\end{equation}
}

The stiffness maximization problem is formulated to maximize the strain energy $W$ of the structure, with the constraint of maximum volume fraction $V_{\text{max}}$, under a predefined displacement boundary condition, as follows:

\begin{equation}
    \begin{aligned}
    \underset{\gamma}{\mathrm{maximize}}\quad & W = \frac{1}{2} \int_\Omega \boldsymbol{\sigma} : \boldsymbol{\varepsilon} \, \mathrm{d}\Omega,\\
    \makebox[2cm][r]{subject to} \quad & V = \frac{1}{\int_\Omega \mathrm{d}\Omega} \int_{\Omega} \tilde{\gamma} \, \mathrm{d}\Omega \ \leq \ V_{\text{max}}, \\
    \quad & \gamma \in [ 0, 1 ],
    \end{aligned}
\label{eq:mono_stiffness}
\end{equation}
\noindent
where $V$ is the volume fraction of the optimized structure.
Note that this $V$ represents the volume fraction of the design variables, which is not necessarily equal to the volume fraction of the actual TPMS structure generated in the de-homogenization process.


\section{Numerical implementation}
\label{sec:impl}

\begin{figure*}[t]
  \centering
  \includegraphics[width=\textwidth]{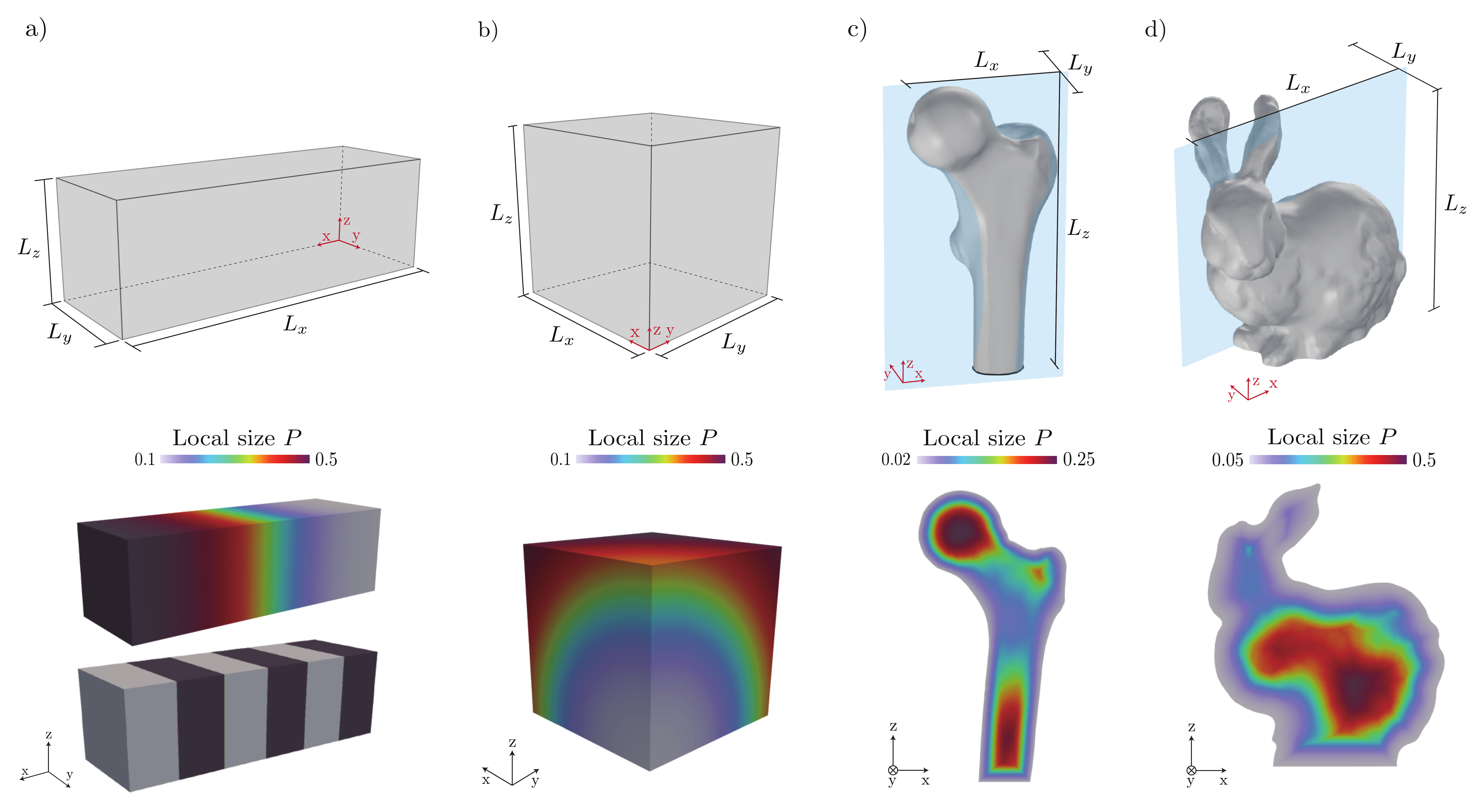}
  \caption{De-homogenization examples: a) 1D-like case with smooth and discrete size distributions, b) simple 3D case, c) complex 3D case with a 3D-scanned bone~\cite{mitsuhashi2008bodyparts3d}, d) complex 3D case with Stanford bunny~\cite{turk1994zippered}. Top row shows the dimensions of the design domain and the positions of the origin with red axes; bottom row shows the corresponding size distributions.}
  \label{fig:examples}
\end{figure*}


COMSOL Multiphysics 6.2 is used for solving the governing equations in the homogenized model and the static structural analysis of the de-homogenized model, as well as for performing the sensitivity analysis for the TO process.
MATLAB R2022b is used for the TO process, including the filtering and the updating of the design variable field, using the COMSOL LiveLink for MATLAB.
Python 3.10.12 is used for the de-homogenization process and the post-processing of the experimental data: the SciPy 1.15.3 is used for computing Gaussian smoothing and solving the Poisson's equation, and the CuMCubes 0.0.3 is used for extracting the solid region of the TPMS structure from the optimized phase distributions.
For the CuMCubes, we accelerate the computation by utilizing GPUs with the PyTorch 2.6.0+cu124, which significantly reduces the computation time.
The numerical implementation is performed on a workstation with an AMD Ryzen Threadripper PRO 7965WX CPU, 512 GB of RAM, and an NVIDIA RTX 4000 Ada GPU.

\subsection{Homogenization}
\label{subsec:homogenization}

As a type of TPMS structure for the demonstration of the proposed method, we used the Gyroid structure, which is commonly used in lattice core design due to its favorable mechanical properties and manufacturability.

In the homogenization process, we use $N_\text{RUC}$ RUCs, where $N_\text{RUC}=10$ in this study, with different sizes $P^{(k)}$ ($k=1,2,\dots,N_\text{RUC}$).
The size $P^{(k)}$ ranges from $P_\text{min} = 5$ [mm] to $P_\text{max} = 20$ [mm], which covers the range of the size distribution optimized by the homogenization-based TO, as follows:
\begin{equation}
  P^{(k)} = P_\text{min} + (P_\text{max}-P_\text{min}) \frac{k-1}{N_\text{RUC}-1} \ (k=1,2,\dots,N_\text{RUC}).
  \label{eq:period_ruc}
\end{equation}

The RUCs are constructed by discretizing the implicit function of the Gyroid structure with a uniform grid of 64 $\times$ 64 $\times$ 64 elements for cell size $P^{(k)}=5,6.7,8.3,10$ [mm] ($k=1,2,3,4$), and 80 $\times$ 80 $\times$ 80 elements for $P^{(k)}=11.7,13.3,15,16.7,18.3,20$ [mm] ($k=5,6,7,8,9,10$).
This discretization is chosen to ensure the sufficient resolution to capture the geometry of the Gyroid structure even at the largest size.
After discretization, the solid region of the Gyroid structure is extracted by calculating the marching cubes algorithm~\cite{lorensen1987marching} based on Eq.~\eqref{eq:tpms_solid}.
The resulting solid region is converted into a quadratic serendipity tetrahedral mesh with a maximum element size of 0.65 [mm] for the FEA.
Then the FEA is performed on the converted mesh to calculate the effective properties of the Gyroid structure for each RUC.

The solid material properties used for the homogenization process are set as follows: Young's modulus $E_\text{solid} = 2450$ [MPa], Poisson's ratio $\nu_\text{solid} = 0.35$, and density $\rho_\text{solid} = 1210$ [kg/m$^3$].
These material properties are defined according to the 3D printing material (Formlabs Clear Resin V5) used for the experimental validation, and thus are also used for the static analysis of the de-homogenized structure as well for consistency.

After performing the homogenization process for each RUC, we obtain the effective Young's modulus $E^{(k)}$, the effective Poisson's ratio $\nu^{(k)}$, and the effective density $\rho^{(k)}$ for each size $P^{(k)}$.
The effective properties of $E^*, \nu^*, \rho^*$ are interpolated as a function of the normalized size $\bar{P}$ using the polynomial fitting technique with the least-squares method, as follows:
\begin{equation}
  \begin{aligned}
    E^* & = a_E \bar{P}^4 + b_E \bar{P}^3 + c_E \bar{P}^2 + d_E \bar{P} + e_E, \\
    \nu^* & = a_\nu \bar{P}^4 + b_\nu \bar{P}^3 + c_\nu \bar{P}^2 + d_\nu \bar{P} + e_\nu, \\
    \rho^* & = a_\rho \bar{P}^4 + b_\rho \bar{P}^3 + c_\rho \bar{P}^2 + d_\rho \bar{P} + e_\rho,
  \end{aligned}
  \label{eq:effective_properties_interpolation}
\end{equation}
\noindent where $a_p, b_p, c_p, d_p, e_p$ are the coefficients of the polynomial fitting for property $p$ ($p=E,\nu,\rho$).
The fitted coefficients are listed in Table~\ref{tab:effective_properties_coefficients}.
Fig.~\ref{fig:interpolation} illustrates the interpolation of the effective properties of the Gyroid structure based on the RUC analysis.
Since the normalized size $\bar{P}$ is defined in the range of $[0,1]$, the interpolation functions can be used to calculate the effective properties for the design variable field $\gamma$ in the TO.


\begin{table}[t]
  \centering
  \caption{Coefficients of the polynomial fitting for the effective properties.}
  \label{tab:effective_properties_coefficients}
  {
    \footnotesize
    \renewcommand{\arraystretch}{1.1}
    \begin{tabular}{c|c|c|ccccc}
    \hline
    Property & Unit & $R^2$ & $a_p$ & $b_p$ & $c_p$ & $d_p$ & $e_p$ \\
    \hline
    $E^*$ & MPa & 1.0 & 624 & -1662 & 1672 & -841 & 263 \\
    $\nu^*$ & - & 1.0 & -0.058 & 0.151 & -0.146 & 0.073 & 0.326 \\
    $\rho^*$ & kg/m$^3$ & 1.0 & 521 & -1442 & 1555 & -887 & 341 \\
    \hline
    \end{tabular}
  }  
\end{table}

\subsection{Topology optimization}
\label{subsec:topology}

Fig.~\ref{fig:simulation_setup} illustrates analysis setups of a numerical example to be considered: the homogenization-based TO on the left, and the performance evaluation on the de-homogenized model on the right.
As shown in the left of Fig.~\ref{fig:simulation_setup}, the design domain is defined as a box with $Lx \times Ly \times Lz$, where $L_x=L_y=L_z=62.5$ [mm] in this study.
Non-design domains are defined as top and bottom layers of the box with a thickness of $Lz/20$, in addition to a small button-like box at the center of the top layer with a size of $Lx/5 \times Ly/5 \times Lz/20$ [mm].
The design domain is meshed into $N_\text{ele} = 10 \times 10 \times 10$ hexahedral elements, while the non-design domain is meshed so as to keep the connectivity between the design domain and the non-design domain, with an element size half to that of the design domain.
On these meshes, the design variable field $\gamma$ is discretized into linear hexahedral elements, whereas the displacement field $\mb{U}$ is discretized into quadratic serendipity hexahedral elements for both non-design domain and design domain.

The bottom face of the analysis domain is fixed, whereas the top face of the button-like box is subjected to a uniform downward displacement of $L_z/10=6.25$ [mm].
This boundary condition is chosen to evaluate the effectiveness of the proposed method in obtaining a high-stiffness structure under locally complex loading conditions, including compression, bending, and shear, as well as to be consistent with the experimental setup described in Section~\ref{sec:exp}.
The identical boundary condition is applied to the de-homogenized model for the performance evaluation, as shown in the right of Fig.~\ref{fig:simulation_setup}.

The solid material properties, $E_\text{solid}, \nu_\text{solid}, \rho_\text{solid}$, used in the homogenization process are assigned to the non-design domain, while the effective properties $E^*, \nu^*, \rho^*$ obtained from the homogenization process are assigned to the design domain, as functions of the design variable field $\gamma$.

The design variable field $\gamma$ is updated iteratively using the Method of Moving Asymptotes (MMA)~\cite{svanberg1987method}, based on the sensitivity of the objective function with respect to the design variable field.
The sensitivity is calculated using the chain rule and the automatic differentiation implemented in COMSOL.
The maximum volume fraction $V_{\text{max}}$ is set to 0.8, and the filter radius $r_{\text{f}}$ is set to $(P_\text{min}+P_\text{max})/2=12.5$ [mm].
The initial design variable field $\gamma$ is set to $0.5$ for all elements in the design domain, which corresponds to the middle size of the TPMS structure with a size of $P=12.5$ [mm].
The maximum number of iterations is set to 100 after confirming the convergence of the optimization process in preliminary tests.

\subsection{De-homogenization}
\label{subsec:dehomogenization}

\subsubsection{Primitive distortion reduction: Gaussian smoothing}

In the conventional PM method, as the most primitive way to reduce the distortion of the de-homogenization, Gaussian smoothing is often applied to the size distribution before applying the PM method, which can be expressed as follows:
\begin{equation}
\tilde{P}^{i,j,k} = \sum_{p=-K}^{K} \sum_{q=-K}^{K} \sum_{r=-K}^{K} W^{p,q,r}\, P^{i+p,j+q,k+r}.
\label{eq:gaussian_filter_discrete}
\end{equation}
Hereinafter, we assume that the design domain is discretized into a regular grid of $N_x \times N_y \times N_z$ points, and the indices $i,j,k$ correspond to the grid points in the $x$-, $y$-, and $z$-directions, respectively.
$W^{p,q,r}$ is the normalized Gaussian kernel, defined as follows:
\begin{equation}
W^{p,q,r} = \frac{ \exp\!\left( -\dfrac{p^2+q^2+r^2}{2\sigma_c^2} \right) }{ \displaystyle \sum_{p'=-K}^{K}\sum_{q'=-K}^{K}\sum_{r'=-K}^{K} \exp\!\left( -\dfrac{{p'}^2+{q'}^2+{r'}^2}{2\sigma_c^2} \right) }.
\label{eq:gaussian_kernel_discrete}
\end{equation}
Here, $\sigma_c$ denotes the standard deviation of the Gaussian kernel measured in the unit of grid cells, and in the present study it was set to
\begin{equation}
\sigma_c = \alpha \max(N_x,N_y,N_z),
\label{eq:sigma_cells}
\end{equation}
where $\alpha$ is a scaling factor that controls the degree of smoothing.
$K$ denotes the truncation radius of the discrete Gaussian kernel in the grid-index space, and the convolution is evaluated over the neighboring cells within $[-K,K]$ in each coordinate direction, with the radius chosen as the smallest integer greater than or equal to $3\sigma_c$.
The boundary values outside the computational domain were handled by nearest-value extension.

\begin{figure}[t]
  \centering
  \includegraphics[width=\columnwidth]{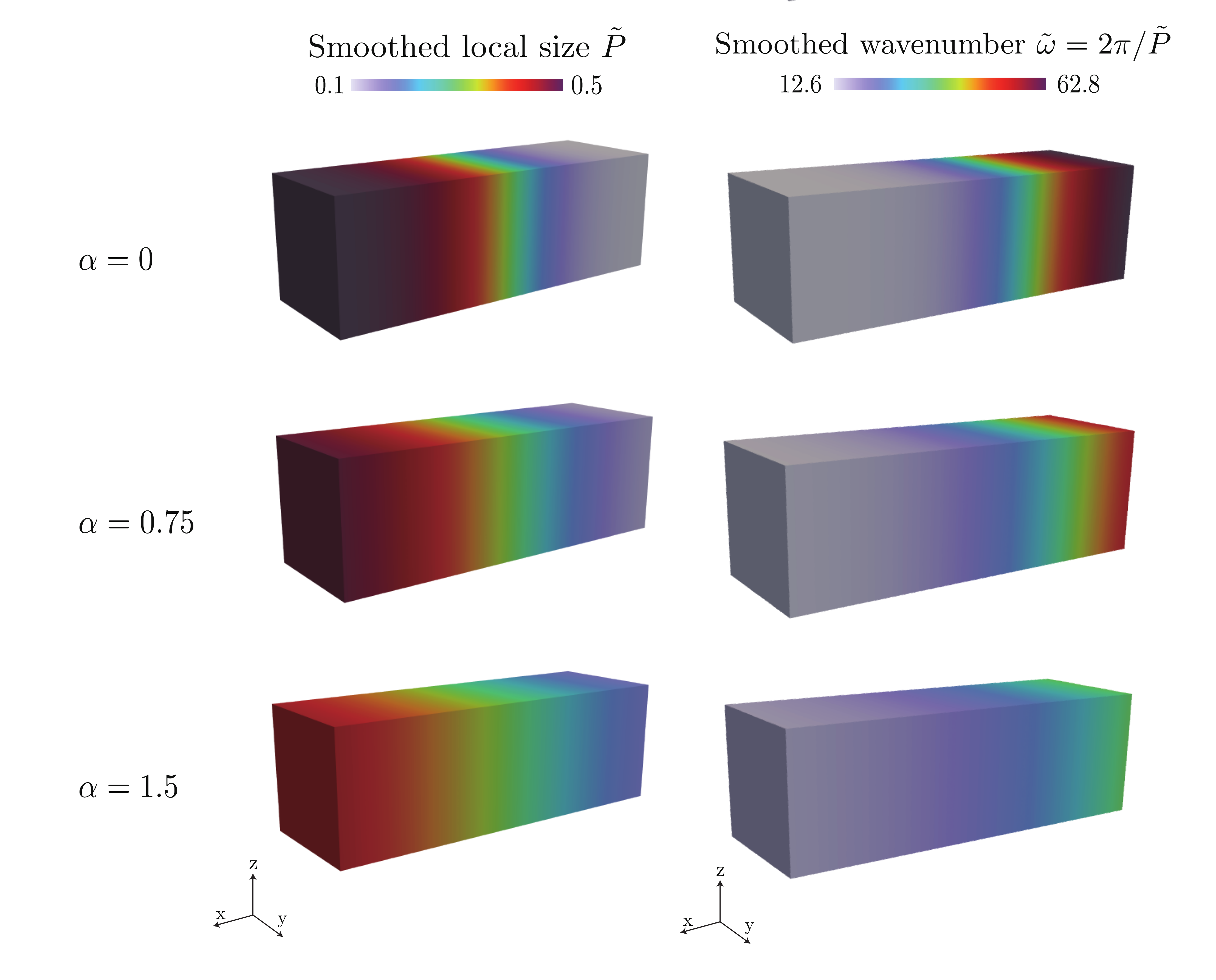}
  \caption{The smoothed size distributions (left) and the corresponding wavenumber distributions (right) for the 1D case. The top row corresponds to the distributions without Gaussian smoothing, which are used for the proposed method. The second and third rows correspond to the distributions with Gaussian smoothing with $\alpha=0.75$ and $\alpha=1.5$, respectively.}
  \label{fig:test1D_size_waveno}
\end{figure}

\begin{figure}[t]
  \centering
  \includegraphics[width=\columnwidth]{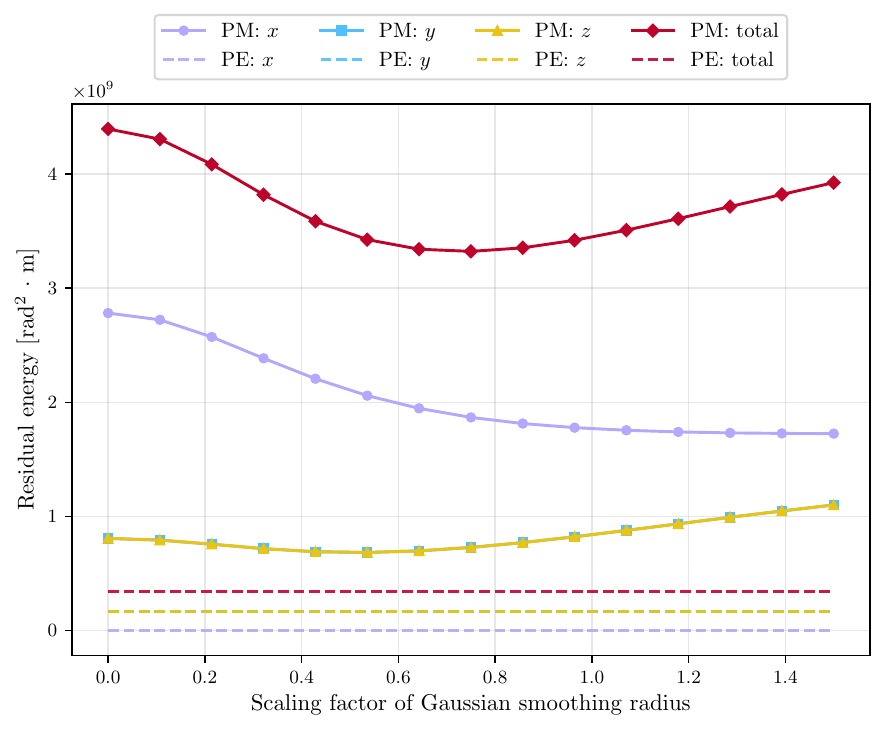}
  \caption{Residual energy vs. scaling factor of Gaussian smoothing $\alpha$ from the conventional periodic modulation (PM) method and the proposed method solving Poisson equation (PE) for $x$-, $y$-, and $z$-directions and the total. For the PE method, the size distribution is not smoothed, and thus the results are plotted as a constant line at all values of $\alpha$.}
  \label{fig:residual_energy_vs_alpha_r}
\end{figure}

\subsubsection{Efficient solution using discrete cosine transform}
\label{subsubsec:dct}

The key idea of the proposed de-homogenization method is to optimize the phase distributions $\phi_s$ by solving the Poisson's equation with Neumann boundary conditions, as given in Eq.~\eqref{eq:euler_lagrange} separately for each direction $s=x,y,z$.
The Poisson's equation can be solved stably and efficiently using discrete cosine transform (DCT)~\cite{ghiglia1994robust}.

Eq.~\eqref{eq:euler_lagrange} can be written in a discretized form, as follows:

\begin{equation}
    \Delta_\text{d} \phi_s^{i,j,k} = b_s^{i,j,k},
    \label{eq:poisson_regularized}
\end{equation}
\noindent where $\Delta_\text{d}$ is the discrete Laplacian operator.
$b_s^{i,j,k}$ is the right-hand side of the equation that includes the discrete divergence of $\boldsymbol{\omega}_s^{i,j,k}$ from the source term in $\Omega$ and the flux term from the Neumann boundary conditions on $\partial \Omega$, defined as follows:

\begin{equation}
  \begin{aligned}
    b_s^{i,j,k} & = \nabla_\text{d} \cdot \boldsymbol{\omega}_s^{i,j,k} - b_\text{flux}^{i,j,k}, \\
    b_\text{flux}^{i,j,k}
    & =
    \frac{g_x^-}{h_x}\,\chi_{i=0}
    - \frac{g_x^+}{h_x}\,\chi_{i=N_x-1} \\
    & + \frac{g_y^-}{h_y}\,\chi_{j=0}
    - \frac{g_y^+}{h_y}\,\chi_{j=N_y-1} \\
    & + \frac{g_z^-}{h_z}\,\chi_{k=0}
    - \frac{g_z^+}{h_z}\,\chi_{k=N_z-1},
  \end{aligned}
  \label{eq:poisson_rhs}
\end{equation}
\noindent where $h_x, h_y, h_z$ are the grid spacings in the $x$-, $y$-, and $z$-directions, respectively;
$\nabla_\text{d}$ is the discrete gradient operator;
$\chi_{i=0}$ is the indicator function that equals 1 at the boundary where $i=0$ and 0 elsewhere, and similarly for $\chi_{i=N_x-1}$, $\chi_{j=0}$, $\chi_{j=N_y-1}$, $\chi_{k=0}$, and $\chi_{k=N_z-1}$.
$g_x^\pm$, $g_y^\pm$, and $g_z^\pm$ are the normal derivatives of $\phi_s$ on the corresponding boundaries, which can be calculated with the outward unit normal vector $\mathbf{n}$ on each boundary, as follows:

\begin{equation}
  g_x^\pm = \boldsymbol{\omega}_s^{i,j,k} \cdot \mathbf{n}\big|_{\Gamma_x^\pm},\quad
  g_y^\pm = \boldsymbol{\omega}_s^{i,j,k} \cdot \mathbf{n}\big|_{\Gamma_y^\pm},\quad
  g_z^\pm = \boldsymbol{\omega}_s^{i,j,k} \cdot \mathbf{n}\big|_{\Gamma_z^\pm}.
  \label{eq:flux_terms}
\end{equation}

Using the type-II forward DCT, the right-hand side $b_s^{i,j,k}$ can be transformed into the frequency domain, as follows:
\begin{equation}
  \begin{aligned}
    \hat{b}_s^{i,j,k} & = \mathcal{F}\{b_s^{i,j,k}\}, \\
    & = \sum_{i=0}^{N_x-1} \sum_{j=0}^{N_y-1} \sum_{k=0}^{N_z-1} \beta_x (l) \beta_y (m) \beta_z (n) b_s^{i,j,k} X_\text{cos} Y_\text{cos} Z_\text{cos}, \\
    & 0 \leq l \leq N_x-1, \quad 0 \leq m \leq N_y-1, \quad 0 \leq n \leq N_z-1,
  \end{aligned}  
  \label{eq:forward_dct}
\end{equation}

\noindent where subscripts $l,m,n$ denote the frequency domain corresponding to the spatial domain indices $i,j,k$; $\hat{b}_s^{i,j,k}$ is the DCT of $b_s^{i,j,k}$, and $\mathcal{F}\{\cdot\}$ denotes the DCT operation; $X_\text{cos}$, $Y_\text{cos}$, and $Z_\text{cos}$ are the DCT basis functions, defined as follows:
\begin{equation}
  \begin{aligned}
    X_\text{cos} & = 2 \cos \left[ \frac{\pi l}{2N_x} \left( 2i + 1 \right) \right], \\
    Y_\text{cos} & = 2 \cos \left[ \frac{\pi m}{2N_y} \left( 2j + 1 \right) \right], \\
    Z_\text{cos} & = 2 \cos \left[ \frac{\pi n}{2N_z} \left( 2k + 1 \right) \right].
  \end{aligned}
\end{equation}
\noindent
$\beta_s (p)$ is the orthogonal normalization factor for the DCT, as defined below:

\begin{equation}
  \begin{aligned}
    \beta_s (p) & = \begin{cases}
      \frac{1}{\sqrt{N_s}} & \text{if } p=0, \\
      \sqrt{\frac{2}{N_s}} & \text{if } p\geq1,
    \end{cases} \quad (s,p)=(x,l),(y,m),(z,n).
  \end{aligned}  
  \label{eq:dct_normalization}
\end{equation}

The DCT of the phase distribution $\hat{\phi}_s^{i,j,k}$ can be obtained by solving the following equation in the frequency domain:

\begin{equation}
  \begin{aligned}
    & \hat{\phi}_s^{i,j,k} = \frac{\hat{b}_s^{i,j,k}}{\lambda_{i,j,k}},\\
    & \lambda_{i,j,k} = \lambda_i + \lambda_j + \lambda_k, \\
    & \lambda_q = \frac{2-2\cos\left(\frac{\pi q}{N_q}\right)}{h_q^2}, \quad q=i,j,k
  \end{aligned}  
  \label{eq:poisson_solution_dct}
\end{equation}
\noindent where $\lambda_i, \lambda_j, \lambda_k$ are the eigenvalues of the Laplacian operator in the frequency domain, and $\lambda_{i,j,k}$ is the combined eigenvalue for the 3D case.

Finally, the optimized phase distribution $\phi_s^{i,j,k}$ can be obtained by applying the inverse 3D DCT type-II to $\hat{\phi}_s^{i,j,k}$, as follows:

\begin{equation}
  \begin{aligned}
    \phi_s^{i,j,k} & = \mathcal{F}^{-1}\{\hat{\phi}_s^{i,j,k}\}, \\
    & = \frac{1}{N_x N_y N_z} \sum_{l=0}^{N_x-1} \sum_{m=0}^{N_y-1} \sum_{n=0}^{N_z-1} \beta_x (l) \beta_y (m) \beta_z (n) \hat{\phi}_s^{l,m,n} X_\text{cos} Y_\text{cos} Z_\text{cos}, \\    
  \end{aligned}
  \label{eq:inverse_dct}
\end{equation}
\noindent where $\mathcal{F}^{-1}\{\cdot\}$ denotes the inverse DCT operation.

After obtaining the optimized phase distributions $\phi_s$, the solid region of the TPMS structure can be extracted by applying the marching cubes algorithm based on Eq.~\ref{eq:tpms_solid}, similar to the process for RUCs in the homogenization step.

\begin{figure*}[t]
  \centering
  \includegraphics[width=\textwidth]{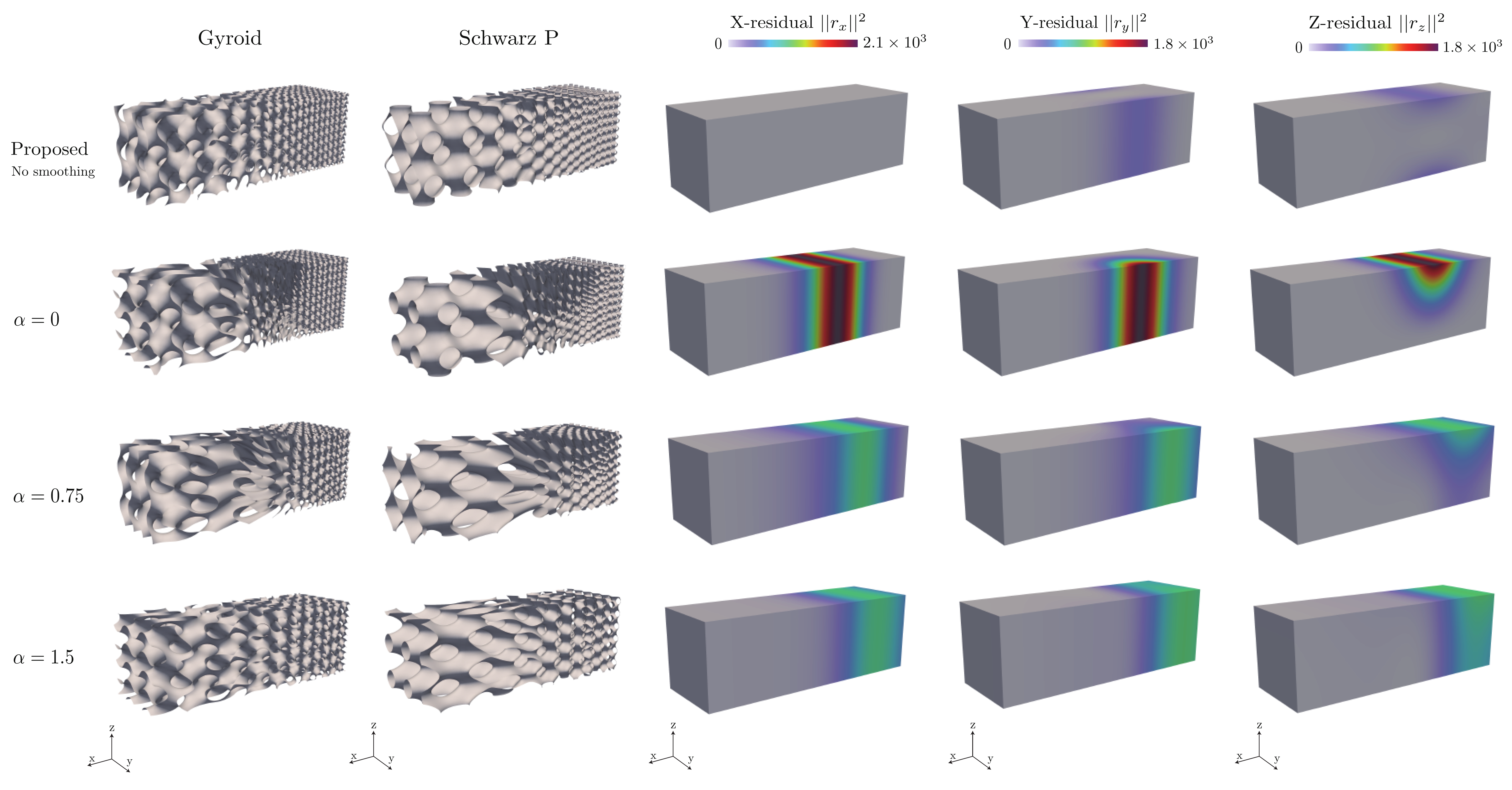}
  \caption{De-homogenized structures and residual distributions for the 1D case: proposed method (first row), PM method with Gaussian smoothing with $\alpha=0$ (second row), $\alpha=0.75$ (third row), and $\alpha=1.5$ (fourth row).}
  \label{fig:test1D_dehomo_residuals}
\end{figure*}

\begin{figure}[t]
  \centering
  \includegraphics[width=\columnwidth]{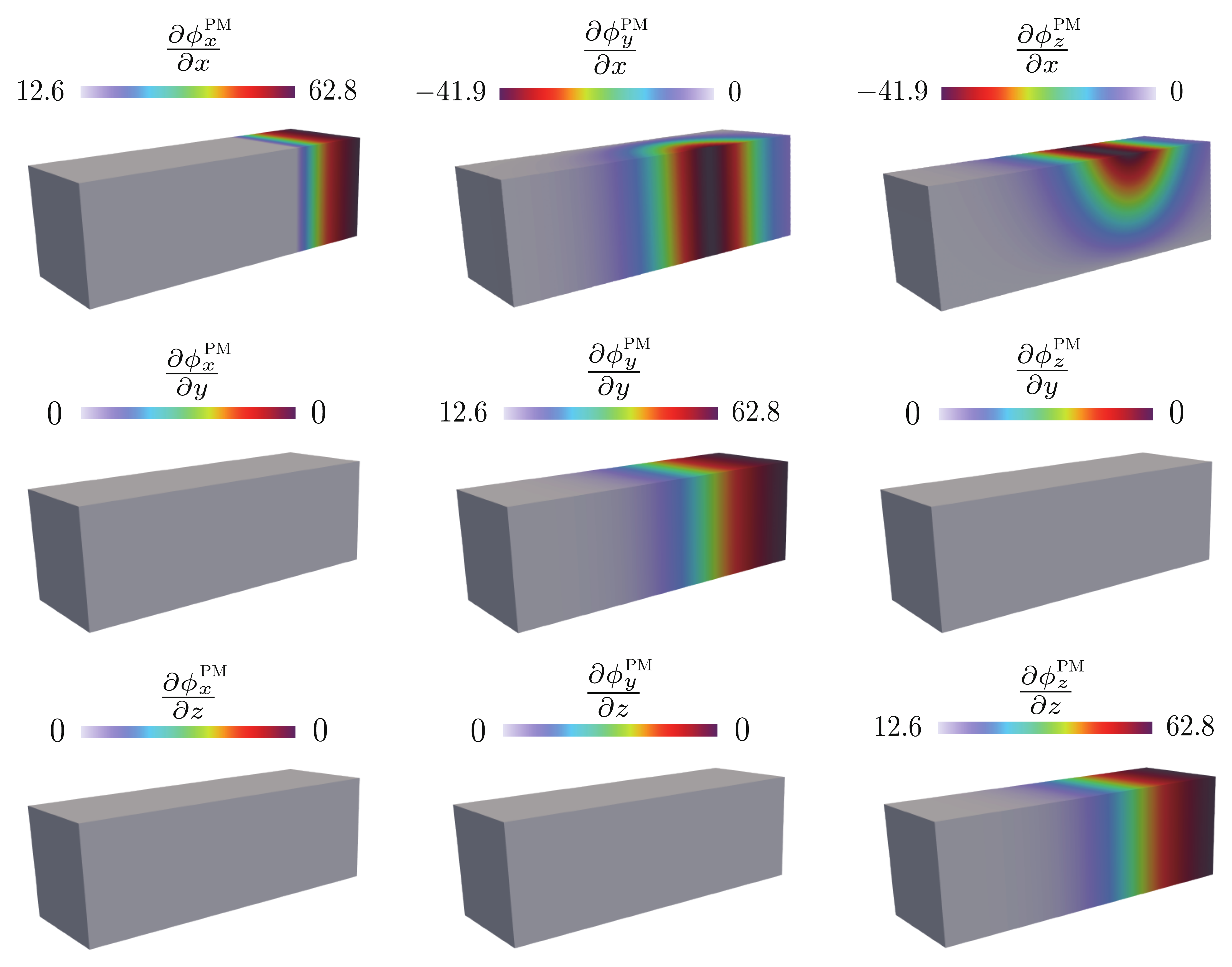}
  \caption{The distributions of the gradients of the phase distributions, showing each component of $\nabla\phi_x^\text{PM}$, $\nabla\phi_y^\text{PM}$, and $\nabla\phi_z^\text{PM}$ for the 1D case from the PM method.}
  \label{fig:test1D_grad_phases}
\end{figure}

\begin{figure}[t]
  \centering
  \includegraphics[width=\columnwidth]{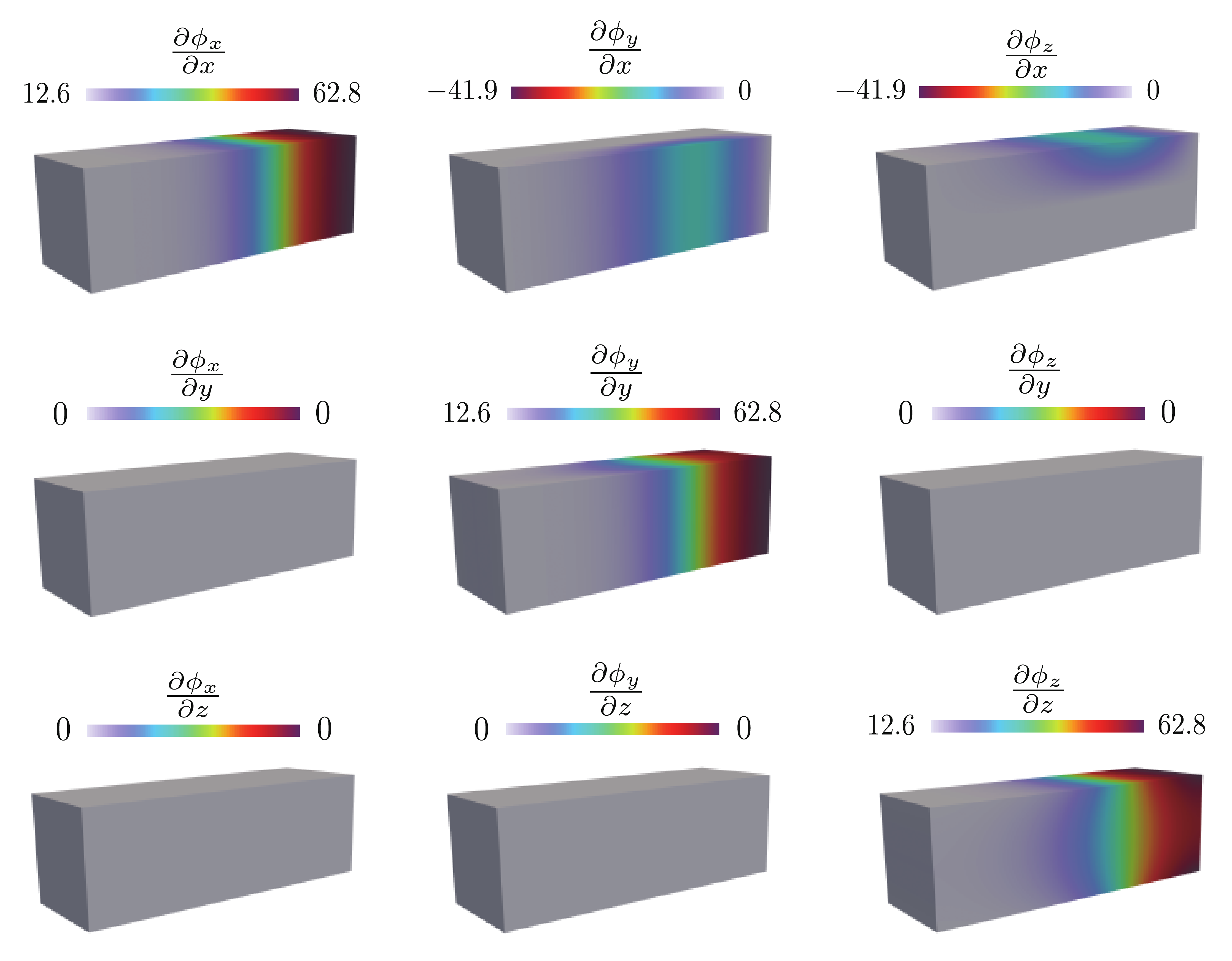}
  \caption{The distributions of the gradients of the phase distributions, showing each component of $\nabla\phi_x$, $\nabla\phi_y$, and $\nabla\phi_z$ for the 1D case from the proposed method.}
  \label{fig:test1D_grad_phases_opt}
\end{figure}

\begin{figure*}[t]
  \centering
  \includegraphics[width=\textwidth]{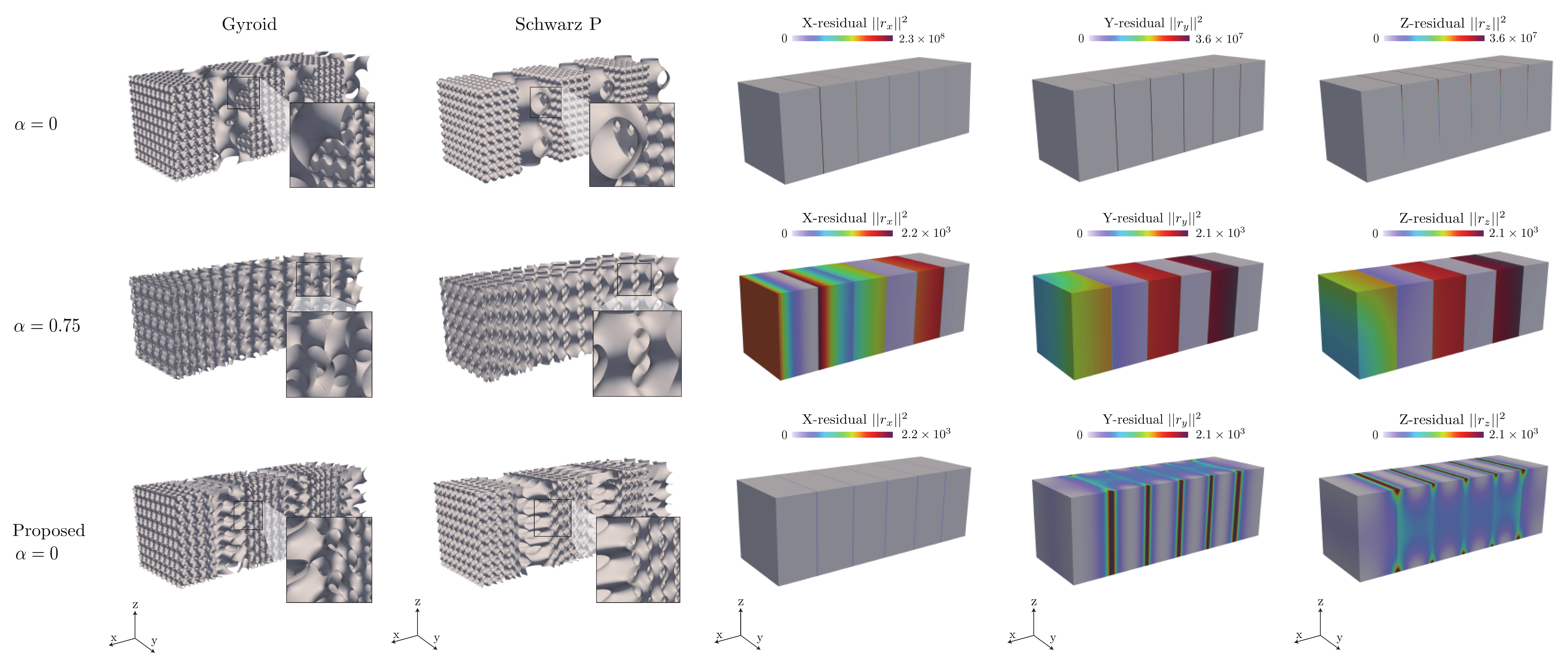}
  \caption{Dehomogenized structures and residual distributions for the 1D case with discrete size distribution: conventional PM method without Gaussian smoothing (top), with Gaussian smoothing with $\alpha=0.75$ (middle), and proposed method (bottom).}
  \label{fig:test1D_discont_dehomo_residuals}
\end{figure*}

\begin{figure*}[t]
  \centering
  \includegraphics[width=\textwidth]{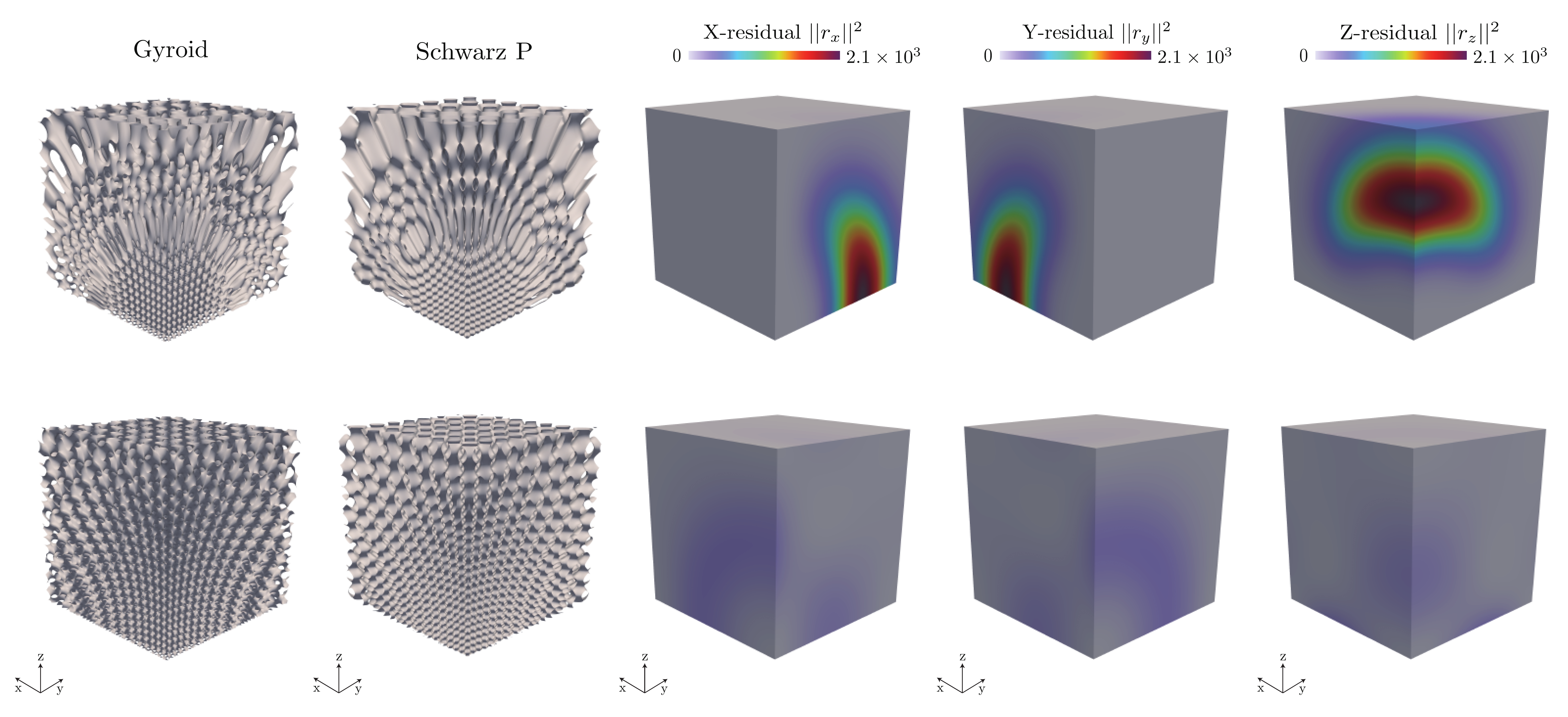}
  \caption{Dehomogenized structures and residual distributions for the 3D case: conventional PM method without Gaussian smoothing (top) and proposed method (bottom).}
  \label{fig:test3D_dehomo_residuals}
\end{figure*}

\subsubsection{De-homogenization examples}
\label{subsubsec:dehomo_examples}

To demonstrate the effectiveness of the proposed de-homogenization method, we apply it to several examples: a simple 1D-like case, a simple 3D case, and two complex 3D cases with 3D-scanned geometries, as shown in Fig.~\ref{fig:examples}.
For each case, the size distribution is defined as a function of the normalized distance field $d^{i,j,k}$, with a sigmoid function as follows:
\begin{equation}
  P^{i,j,k} = P_\text{min} + \frac{P_\text{max}-P_\text{min}}{1+ \exp \left(-\kappa \left(d^{i,j,k}-0.5\right)\right)}.
  \label{eq:size_distribution_1D}
\end{equation}
\noindent where $\kappa$ is a parameter that controls the steepness of the sigmoid function.

The 1D-like case is designed to have a size distribution that varies only in the x-direction, where $d^{i,j,k}$ represents the normalized distance from the origin in the x-direction, i.e., $d^{i,j,k} = x^i / L_x$.
This size distribution is chosen to have a gradual change from $P_\text{min}$ to $P_\text{max}$ along the $x$-direction, but at the same time to have a more intense size distribution at the middle of the design domain, which allows us to evaluate the performance of the proposed method in handling size distributions with varying intensity in one direction.
To demonstrate the robustness of the proposed method against the discrete size distribution, we also apply the proposed method to a discrete size distribution that is obtained with a normalized distance field $d^{i,j,k}$ defined by the following function:
\begin{equation}
  d^{i,j,k} = \begin{cases}
    1.0 & \text{if } n\frac{L_x}{6} \leq x^i < (n+1)\frac{L_x}{6}, \quad (n=0,2,4) \\
    0.0 & \text{else },
  \end{cases}
\end{equation}
\noindent which creates a discrete size distribution with six distinct size regions along the x-direction, with the size distribution alternating between $P_\text{min}$ and $P_\text{max}$ every $L_x/6$.

Similar to the smooth 1D-like case, the 3D case is designed to have a size distribution that varies in all three directions depending on the distance from the origin, defined by substituting the distance function $d^{i,j,k}$ in Eq.~\eqref{eq:size_distribution_1D} with the following function that calculates the normalized distance from the origin in 3D space:
\begin{equation}
    d^{i,j,k} = \sqrt{\frac{(x^i/L_x)^2+(y^j/L_y)^2+(z^k/L_z)^2}{3}}.
  \label{eq:size_distribution_3D}
\end{equation}

In the two complex 3D cases, the 3D-scanned geometries are used as a zeroset surface for the size distribution for each case: a bone STL from BodyParts3D dataset~\cite{mitsuhashi2008bodyparts3d}, and and a bunny STL from Stanford 3D Scanning Repository~\cite{turk1994zippered}.
The size distribution is defined by substituting the distance function $d^{i,j,k}$ in Eq.~\eqref{eq:size_distribution_1D} with the distance field $d_\text{surf}^{i,j,k}$ from the surface of the 3D-scanned geometries.
The parameters used for the size distributions in these examples are listed in Table~\ref{tab:dehomo_parameters}.

\begin{table}[t]
  \centering
  \caption{Parameters used in the de-homogenization examples.}
  \label{tab:dehomo_parameters}
  {
    \footnotesize
    \renewcommand{\arraystretch}{1.1}
    \begin{tabular}{ccccc}
    \hline
    Parameter & 1D-like & 3D & Bone & Bunny \\
    \hline
    $N_x$ & 360 & 360 & 678 & 792 \\
    $N_y$ & 120 & 360 & 457 & 659 \\
    $N_z$ & 120 & 360 & 1132 & 793 \\
    $L_x$ & 3.0 & 3.0 & 1.8 & 3.0 \\
    $L_y$ & 1.0 & 3.0 & 1.2 & 2.5 \\
    $L_z$ & 1.0 & 3.0 & 3.0 & 3.0 \\
    $P_\text{min}$ & 0.1 & 0.1 & 0.02 & 0.05 \\
    $P_\text{max}$ & 0.5 & 0.5 & 0.25 & 0.5 \\
    $\kappa$ & 10.0 & 10.0 & 5.0 & 8.0 \\
    \hline
    \end{tabular}
  }  
\end{table}

\begin{figure*}[t]
  \centering
  \includegraphics[width=\textwidth]{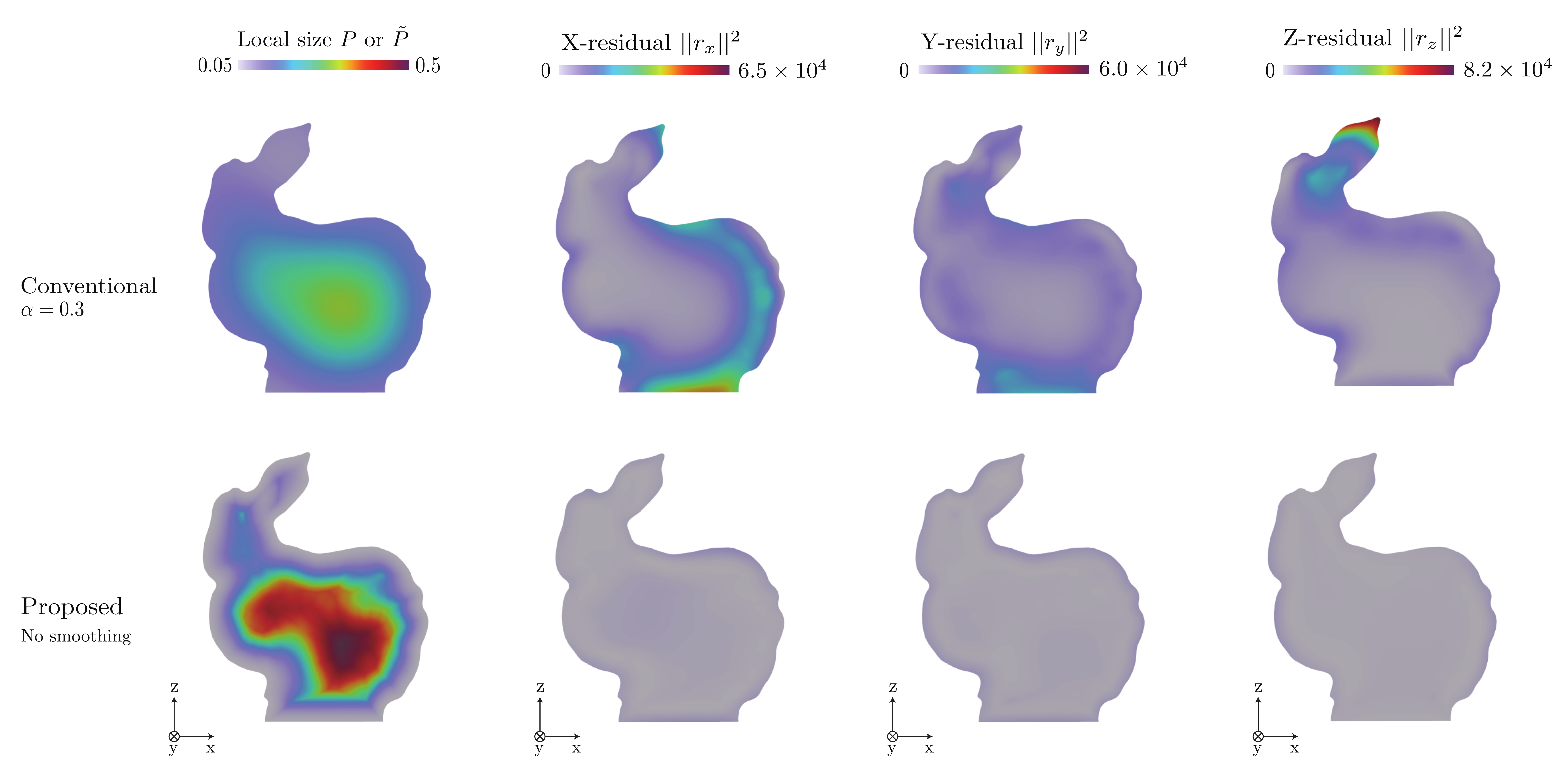}
  \caption{Size distributions and residual distributions for the bunny case: conventional PM method with Gaussian smoothing with $\alpha=0.3$ (top) and proposed method (bottom). The local size distributions correspond to the smoothed size $\tilde{P}$ (top) and the original size $P$ (bottom).}
  \label{fig:bunny_size_residuals}
\end{figure*}

\begin{figure}[t]
  \centering
  \includegraphics[width=\columnwidth]{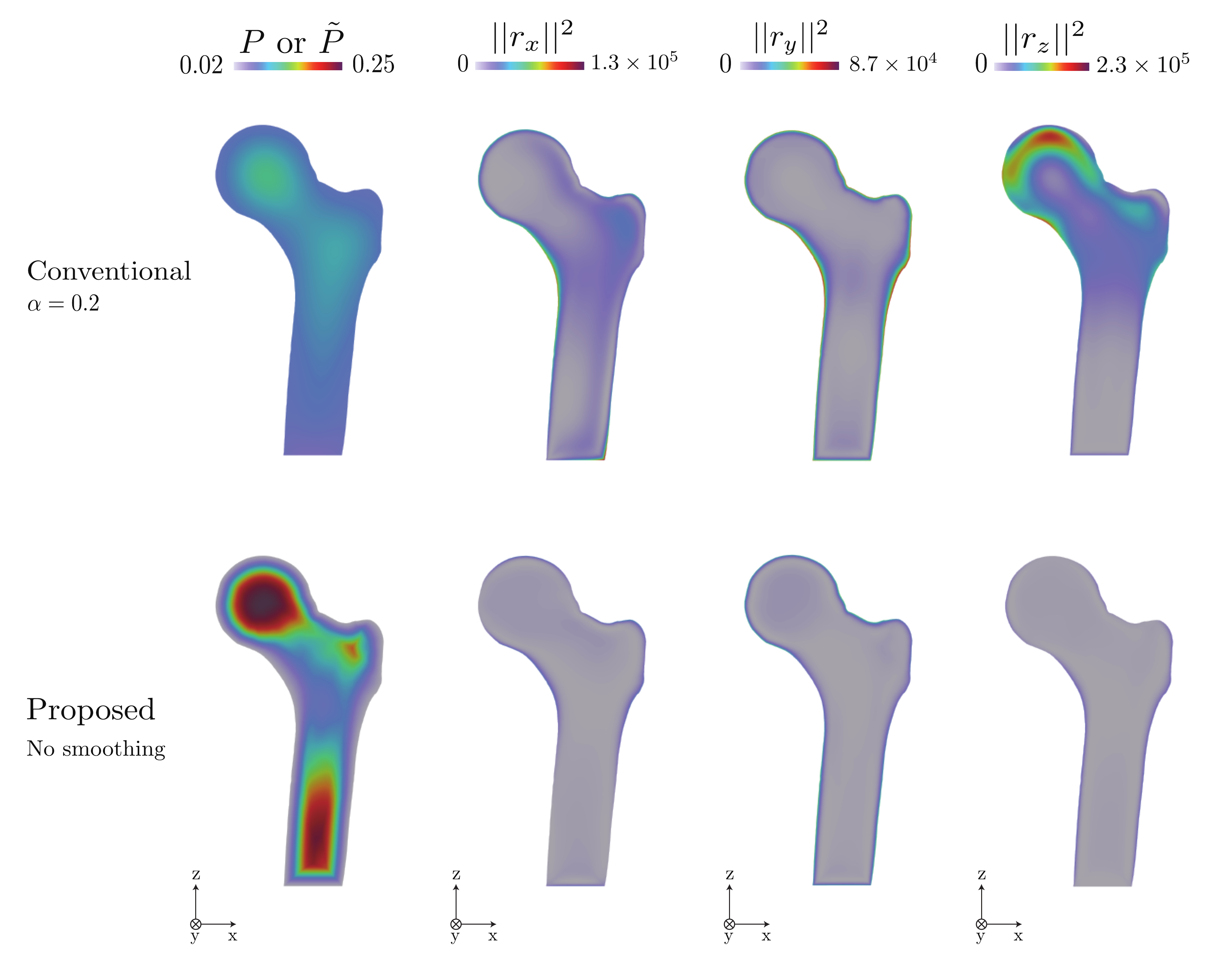}
  \caption{Size distributions and residual distributions for the bone case: conventional PM method with Gaussian smoothing with $\alpha=0.2$ (top) and proposed method (bottom). The local size distributions correspond to the smoothed size $\tilde{P}$ (top) and the original size $P$ (bottom).}
  \label{fig:bone_size_residuals}
\end{figure}

Fig.~\ref{fig:test1D_size_waveno} shows how the Gaussian smoothing with different values of $\alpha$ affects the size distribution and the corresponding wavenumber fields for the 1D-like case.
As the scaling factor $\alpha$ increases, the size distribution becomes smoother, and thus the wavenumber fields become more ambiguous, which can lead to smoother phase distributions from the PM method, since they are simply obtained by multiplying the target wavenumber fields by the spatial coordinates, i.e., $\phi_s=\omega s$.

Additionally, to evaluate the effect of the scaling factor of Gaussian smoothing $\alpha$ on the de-homogenization results from the conventional PM method, we apply the PM method with different values of $0.0 \leq \alpha \leq 1.5$ to the same size distribution, and compare the results with those from the proposed method, as shown in Fig.~\ref{fig:residual_energy_vs_alpha_r}.
The residual energy $||r_s||^2$ is calculated as the sum of the squared differences between the gradients of the optimized phase distributions and the desired wavenumber fields $\omega_s$, or between the modulated phase distributions and the smoothed wavenumber fields $\tilde{\omega}_s$.
The residual energy for the proposed method corresponds to the objective function of the de-homogenization optimization problem $\int_\Omega ||\nabla \phi_s - \omega_s||^2 \, d\Omega$, while the residual energy for the PM method is similarly calculated as $\int_\Omega ||\nabla \phi_s - \tilde{\omega}_s||^2 \, d\Omega$.
Fig.~\ref{fig:residual_energy_vs_alpha_r} shows that the residual energy in $x$ direction from the PM method decreases as the scaling factor $\alpha$ increases, whereas the residual energy in $y$ and $z$ directions stays or increases as $\alpha$ increases, with the total residual energy being minimized at $\alpha=0.75$, resulting in $3.3\times10^9$ rad$^2$.
This indicates that the PM method with a certain level of Gaussian smoothing can achieve a better balance between the residuals in different directions, but the total residual energy is still relatively high.
In contrast, the proposed method achieved a tenfold lower total residual energy, $3.4\times10^8$ rad$^2$, than the minimum total residual energy from the PM method, even without the need for Gaussian smoothing.


Fig.~\ref{fig:test1D_dehomo_residuals} shows the de-homogenization results for the 1D-like case for Gyroid and Schwarz P structures, with the corresponding residual distributions in $x$, $y$, and $z$ directions.
The results show that the proposed method can achieve much lower residuals in all three directions than the PM method with Gaussian smoothing, resulting in much smaller distortions in the de-homogenized structure.
Fig.~\ref{fig:test1D_dehomo_residuals} also shows that the PM method with Gaussian smoothing can reduce the residuals and the distortion around the area with the most intense size change, but instead the de-homogenized structure becomes more stretched entirely because of the smoothed wavenumber fields, occuring more distortion in the areas with less intense size change.

To further analyze the optimized phase distributions, we compare the gradients of the phase distributions from the PM method without Gaussian smoothing and the proposed method, as shown in Fig.~\ref{fig:test1D_grad_phases} and Fig.~\ref{fig:test1D_grad_phases_opt}, respectively.
In the ideal case where the phase distributions perfectly match the target wavenumber field $\omega$ shown in Fig.~\ref{fig:test1D_size_waveno}, the diagonal components of the gradients of the phase distributions, i.e., $\frac{\partial \phi_x}{\partial x}$, $\frac{\partial \phi_y}{\partial y}$, and $\frac{\partial \phi_z}{\partial z}$, match the target wavenumber $\omega$, while the off-diagonal components are zero.
In Fig.~\ref{fig:test1D_grad_phases}, $\frac{\partial \phi_y}{\partial y}$ and $\frac{\partial \phi_z}{\partial z}$ from the PM method closely match the target wavenumber $\omega$, while $\frac{\partial \phi_x}{\partial x}$ shows a significant deviation from the target wavenumber.
At the same time, the $x$-direction off-diagonal components of the gradients of the $y$- and $z$-direction phase distributions, i.e., $\frac{\partial \phi_y}{\partial x}$ and $\frac{\partial \phi_z}{\partial x}$, shows significant non-zero values, which indicates that the phase distributions from the PM method are not properly aligned with the coordinate axes, resulting in distortions in the de-homogenized structure.
In contrast, in Fig.~\ref{fig:test1D_grad_phases_opt}, the gradients of the optimized phase distributions from the proposed method closely match the target wavenumber fields in all three directions, while the off-diagonal components are much closer to zero than those from the PM method, indicating that the optimized phase distributions are well-aligned with the coordinate axes, resulting in much less distortion in the de-homogenized structure.

Fig.~\ref{fig:test1D_discont_dehomo_residuals} shows the de-homogenization results for the 1D-like case with the discrete size distribution, with the corresponding residual distributions in $x$, $y$, and $z$ directions.
The results show that the PM method without Gaussian smoothing completely fails to handle the discrete size distribution, resulting in significant residuals by four or five orders of magnitude higher than the proposed method and the PM method with Gaussian smoothing, and thus discontinuities at the size transition regions in the de-homogenized structure.
On the other hand, the PM method with Gaussian smoothing can reduce the residuals and avoid discontinuities, but it still shows significant residuals in a repetitive pattern, resulting in the de-homogenized structure being stretched entirely losing the intended discrete size distribution.
In contrast, the proposed method can achieve residuals restricted to the narrow bands around the size transition areas, leading to a de-homogenized structure without discontinuities and with a much better reflection of the intended discrete size distribution than the PM method with/without Gaussian smoothing.
However, even in the proposed method, the region with $P_\text{max}$ led to much smaller size in the de-homogenized structure compared to the intended size, i.e., two cells in $y$ and $z$ directions, as observed in Fig.~\ref{fig:test1D_discont_dehomo_residuals}.
Despite the robustness of the proposed method against the discrete size distribution, it is worth noting that this method only can minimize the residuals, and thus too large leaps in the size distribution can still cause the deviation from the intended size distribution in the de-homogenized structure.

Regarding the 3D case, the de-homogenization results from the PM method without Gaussian smoothing and the proposed method are shown in Fig.~\ref{fig:test3D_dehomo_residuals}, with the corresponding residual distributions in $x$, $y$, and $z$ directions.
The residual distributions from the PM method show more significant residuals in the regions that are farther from the origin and have more intense size variations, which is consistent with the theoretical expectation of Eq.~\eqref{eq:grad_phase_spatial}, where the residuals are obtained by multiplying the spatial coordinates by the gradients of the target wavenumber fields.
Similar to the 1D-like case, the results show that the proposed method can achieve much lower residuals in all three directions than the PM method, resulting in much smaller distortions in the de-homogenized structure.
This indicates that the proposed method can effectively handle the size distribution with variations in all three directions, and thus can be applied to a wider range of applications than the PM method with Gaussian smoothing.

For the bunny and bone cases, the PM method strongly requires Gaussian smoothing to avoid discontinuities in the de-homogenized structure because of more complex size distributions with more intense variations.
In other words, if the PM method is applied without Gaussian smoothing, the de-homogenization process requires extremely high resolution to capture the intense variations in the size distribution, which can lead to significant computational cost and memory usage.
In this study, Gaussian smoothing with $\alpha=0.3$ and $\alpha=0.2$ is applied for the bunny and bone cases, respectively.
Fig.~\ref{fig:bunny_size_residuals} and Fig.~\ref{fig:bone_size_residuals} show the size distributions and the residual distributions in $x$, $y$, and $z$ directions on the slice at $y=1.32$ and $y=0.44$ for the bunny and bone cases, respectively.
The slice locations are chosen to include the regions with the most intense size variations, which are around the middle of the design domain in $y$ direction for both cases.

\begin{figure}[t]
  \centering
  \includegraphics[width=\columnwidth]{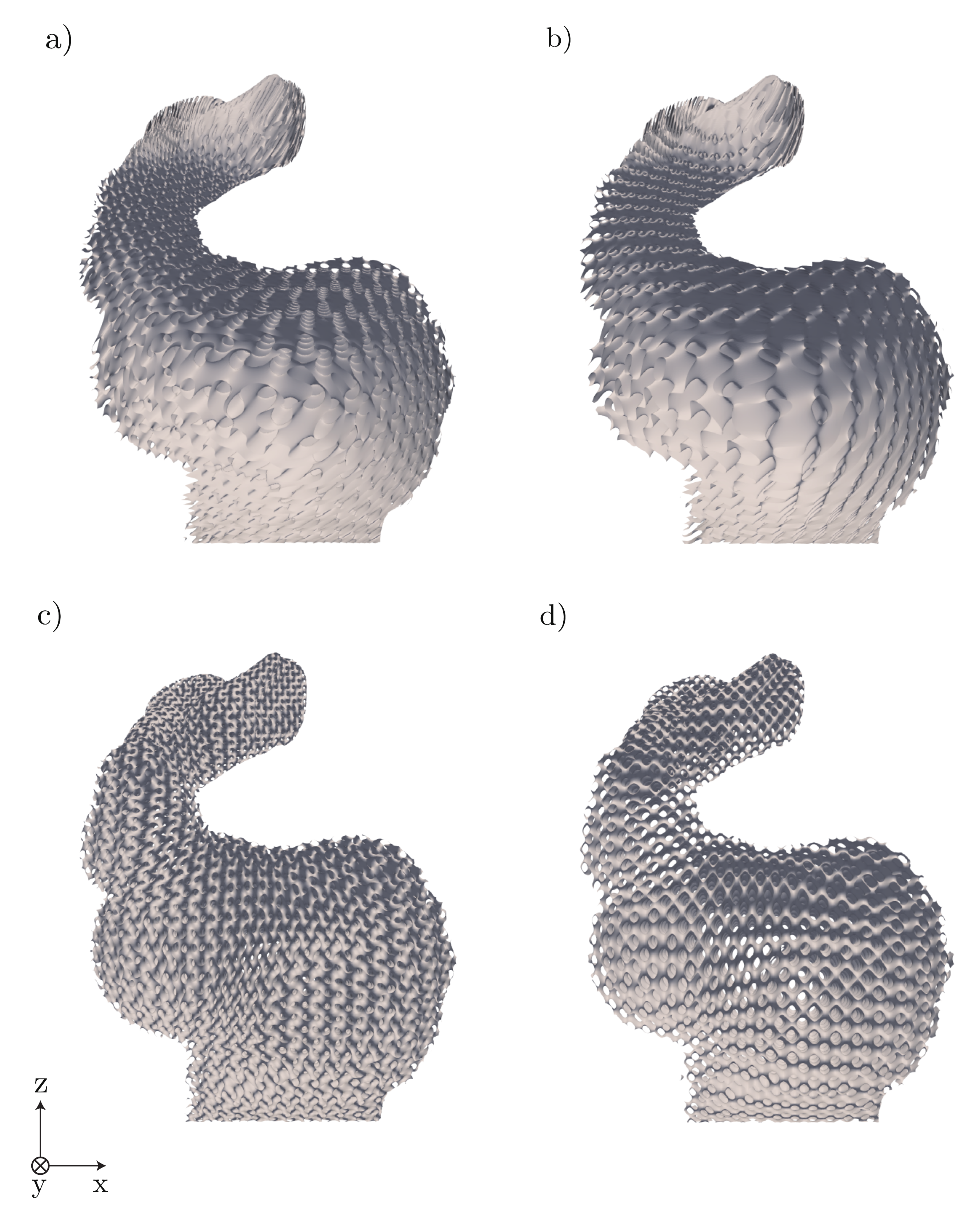}
  \caption{De-homogenized structure for the bunny case: a) Gyroid structure from the conventional PM method, b) Schwarz P structure from the PM method, c) Gyroid structure from the proposed method, d) Schwarz P structure from the proposed method.}
  \label{fig:bunny_dehomo}
\end{figure}

For the bunny case, the smoothed size distribution shows a much more gradual change than the original size distribution, where the regions closer to the center of the bunny body have larger size than the regions farther from the center.
The residual distributions from the PM method exhibit more prominent residuals in the regions farther from the origin, rather than the regions with more intense size variations.
This is because the smooth size distribution after Gaussian smoothing leads to less intense size variations, relatively increasing the effect of the spatial coordinates in the residuals, as expected from Eq.~\eqref{eq:grad_phase_spatial}.
In contrast, the residual distributions from the proposed method show tiny residuals around the bunny skin, indicating that the optimized phase distributions can closely match the target wavenumber fields even in the regions with intense size variations.
Fig.~\ref{fig:bunny_dehomo} shows the de-homogenized structure cut at the same slice as Fig.~\ref{fig:bunny_size_residuals} for the bunny case.
The de-homogenized structure from the PM method is stretched entirely, largely distorting the shape of the TPMS unit cell, especially in the regions farther from the origin with more significant residuals.
On the other hand, the de-homogenized structure from the proposed method can better reflect the target size distribution with much less distortion, where the regions closer to the bunny skin have smaller size than the regions farther from the bunny skin, with the shape of the TPMS unit cell being well preserved.

\begin{figure}[t]
  \centering
  \includegraphics[width=0.9\columnwidth]{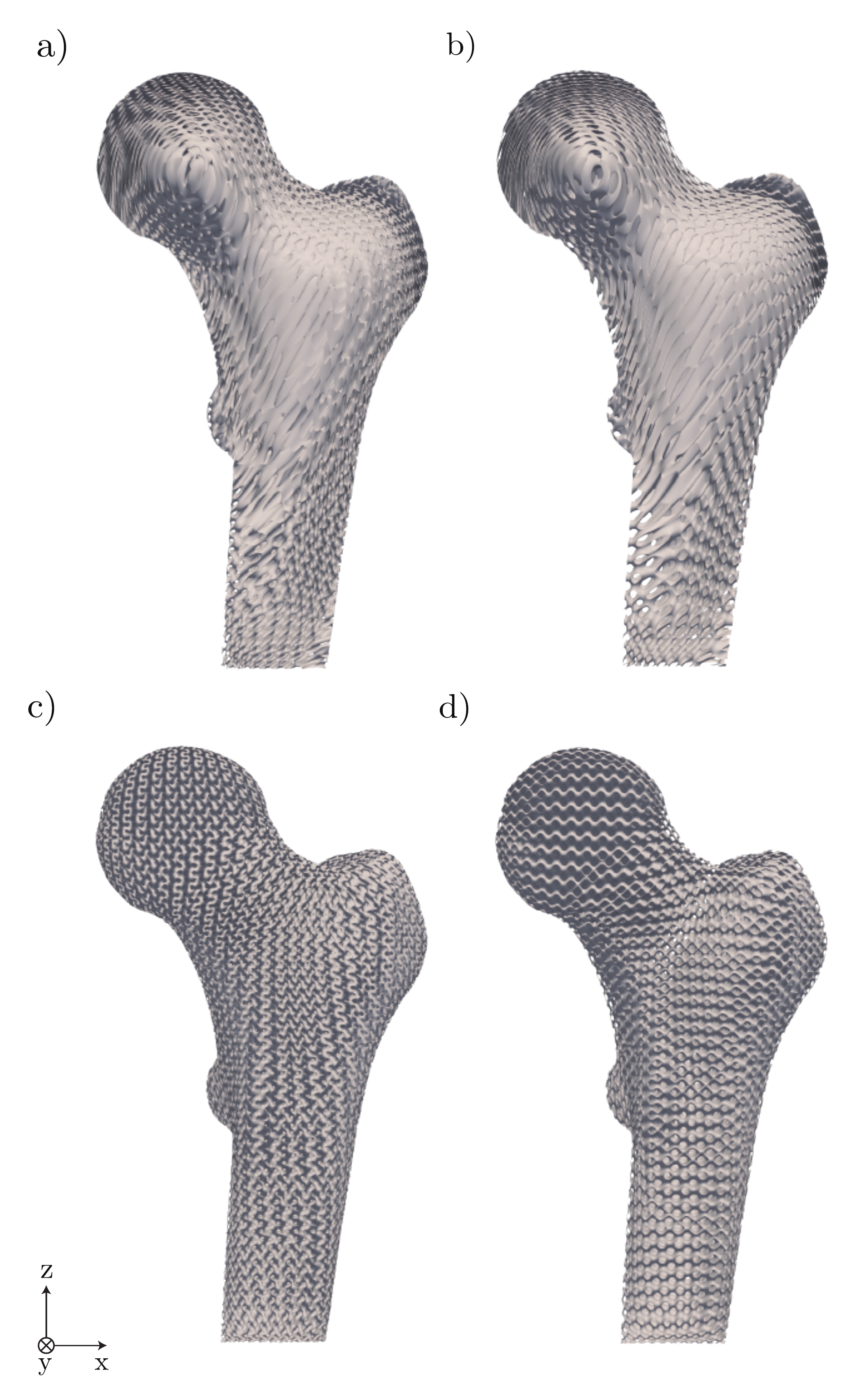}
  \caption{De-homogenized structure for the bone case: a) Gyroid structure from the conventional PM method, b) Schwarz P structure from the PM method, c) Gyroid structure from the proposed method, d) Schwarz P structure from the proposed method.}
  \label{fig:bone_dehomo}
\end{figure}

Similarly, for the bone case, the residual distortions from the PM method are more significant in the regions farther from the origin.
In particular, the top sphere-shaped region, which has more intense size variations and is farther from the origin, exhibits the most significant residuals, leading to a highly distorted de-homogenized structure in that region, as shown in Fig.~\ref{fig:bone_dehomo}.
The de-homogenized structure from the PM method interestingly shows a region with the structure distorted along with the sphere surface, which is unintended and problematic when the structure is optimized assuming the TPMS unit cell shape is preserved.
In contrast, the proposed method successfully achieves much smaller residuals even in the top sphere-shaped region, leading to much less undesired and unintended distortions in the de-homogenized structure, with the unit cell shape being well preserved unlike the PM method.

\section{Numerical examples}
\label{sec:num}

To demonstrate that de-homogenized TPMS structures closely reflect the optimized size distribution, thereby achieving improved structural performance, we applied the proposed method to the lattice core design problem described in Section~\ref{sec:impl}.

The thickness of the solid wall in the de-homogenized structure is determined by a fixed $c$ value in Eq.~\eqref{eq:tpms_solid} for simplicity, defined as follows:
\begin{equation}
  c = \sqrt{2} P \sin \left( \frac{\pi t}{P} \right),
   \label{eq:thickness_specification}
\end{equation}
\noindent where $t$ is the desired thickness of the solid wall at a specific location, set as $t = 0.5$ [mm] in this study, as detailed in~\ref{ap:thickness}.

\begin{figure*}[t]
  \centering
  \includegraphics[width=\textwidth]{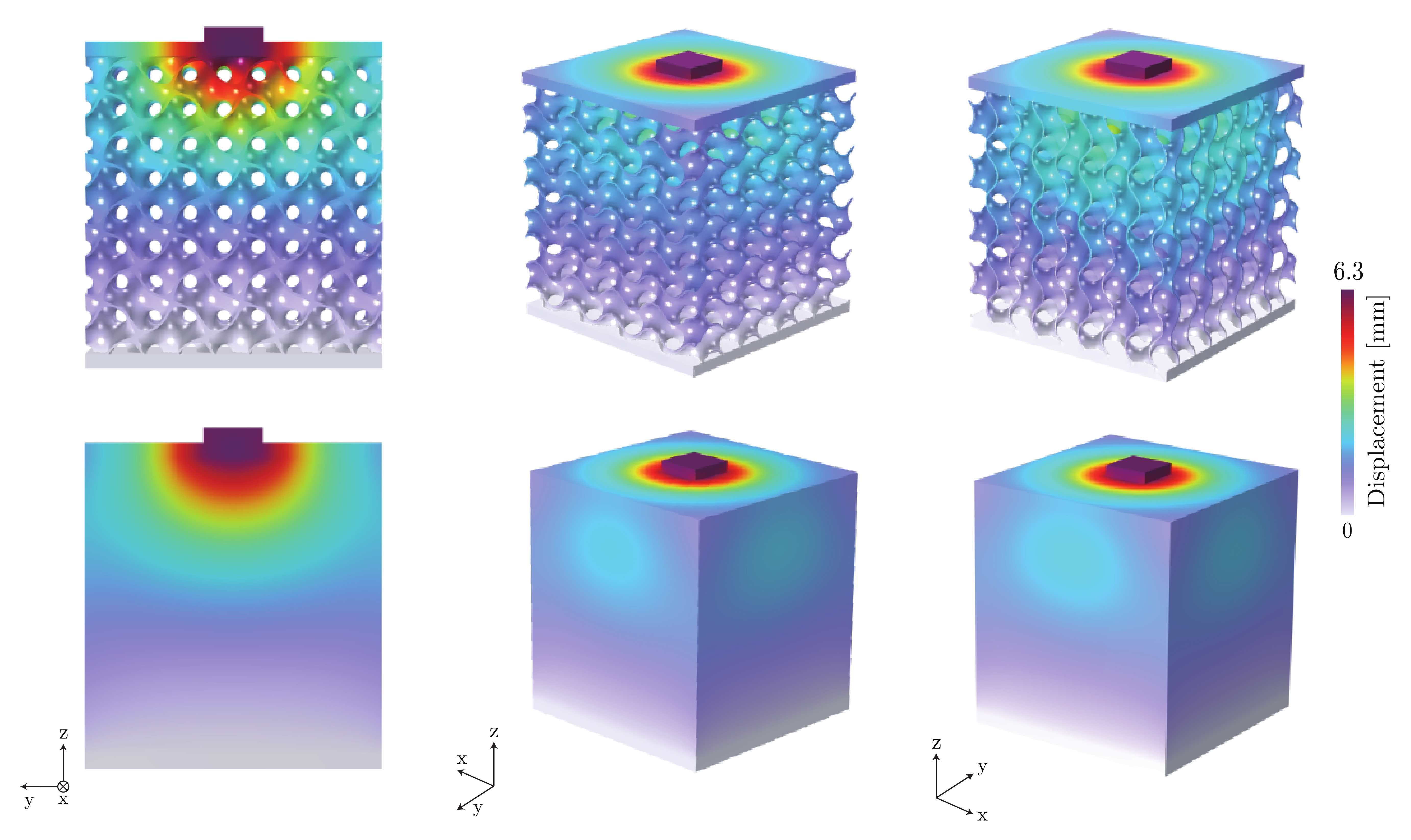}
  \caption{Displacement distribution for the uniform size case with the sectional view at the center and two opposite side views of the de-homogenized model (top) and homogenized model (bottom).}
  \label{fig:uniform_disp}
\end{figure*}

\begin{figure*}[t]
  \centering
  \includegraphics[width=\textwidth]{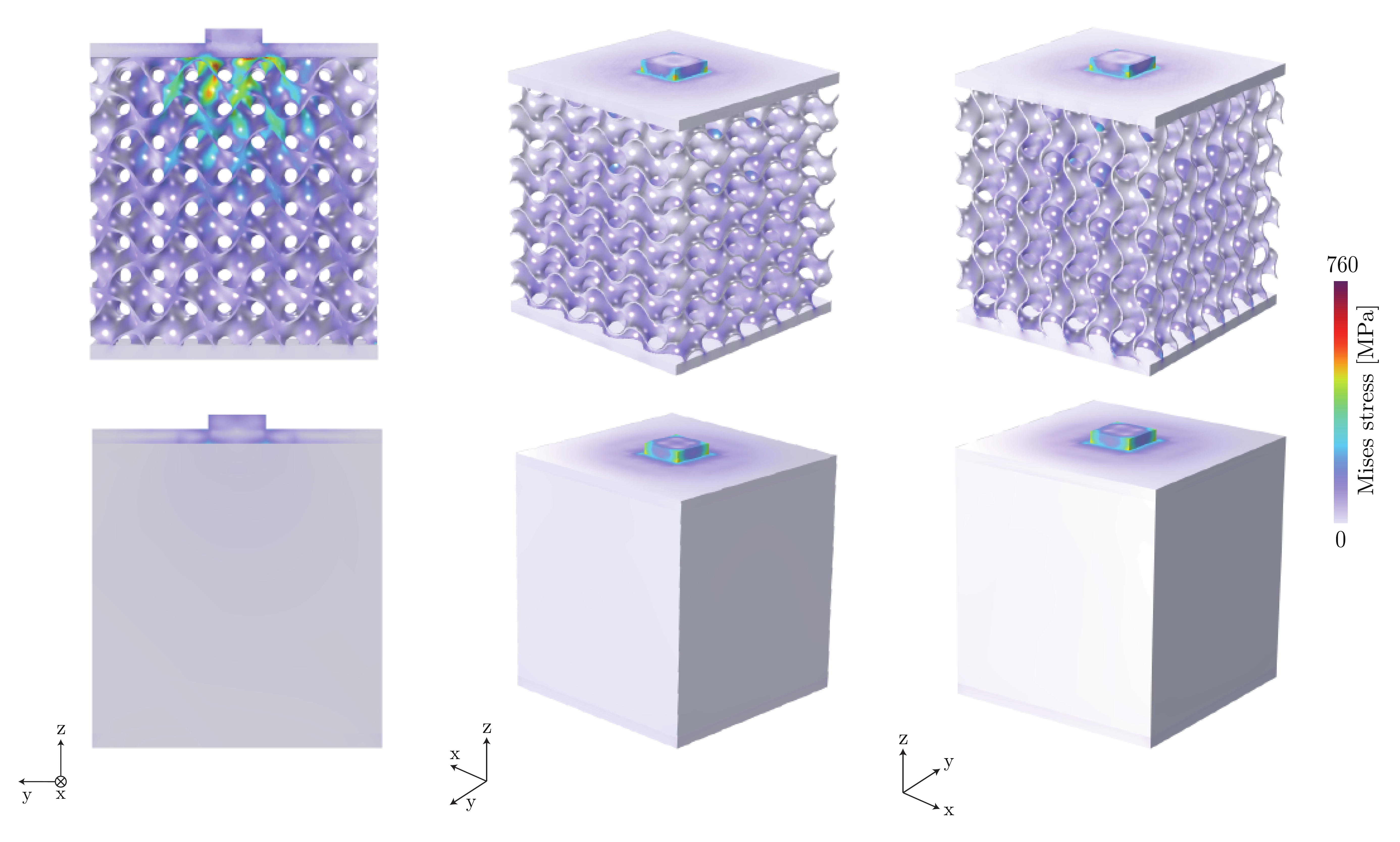}
  \caption{Stress distribution for the uniform size case with the sectional view at the center and two opposite side views of the de-homogenized model (top) and homogenized model (bottom).}
  \label{fig:uniform_stress}
\end{figure*}

\begin{figure}[t]
  \centering
  \includegraphics[width=\columnwidth]{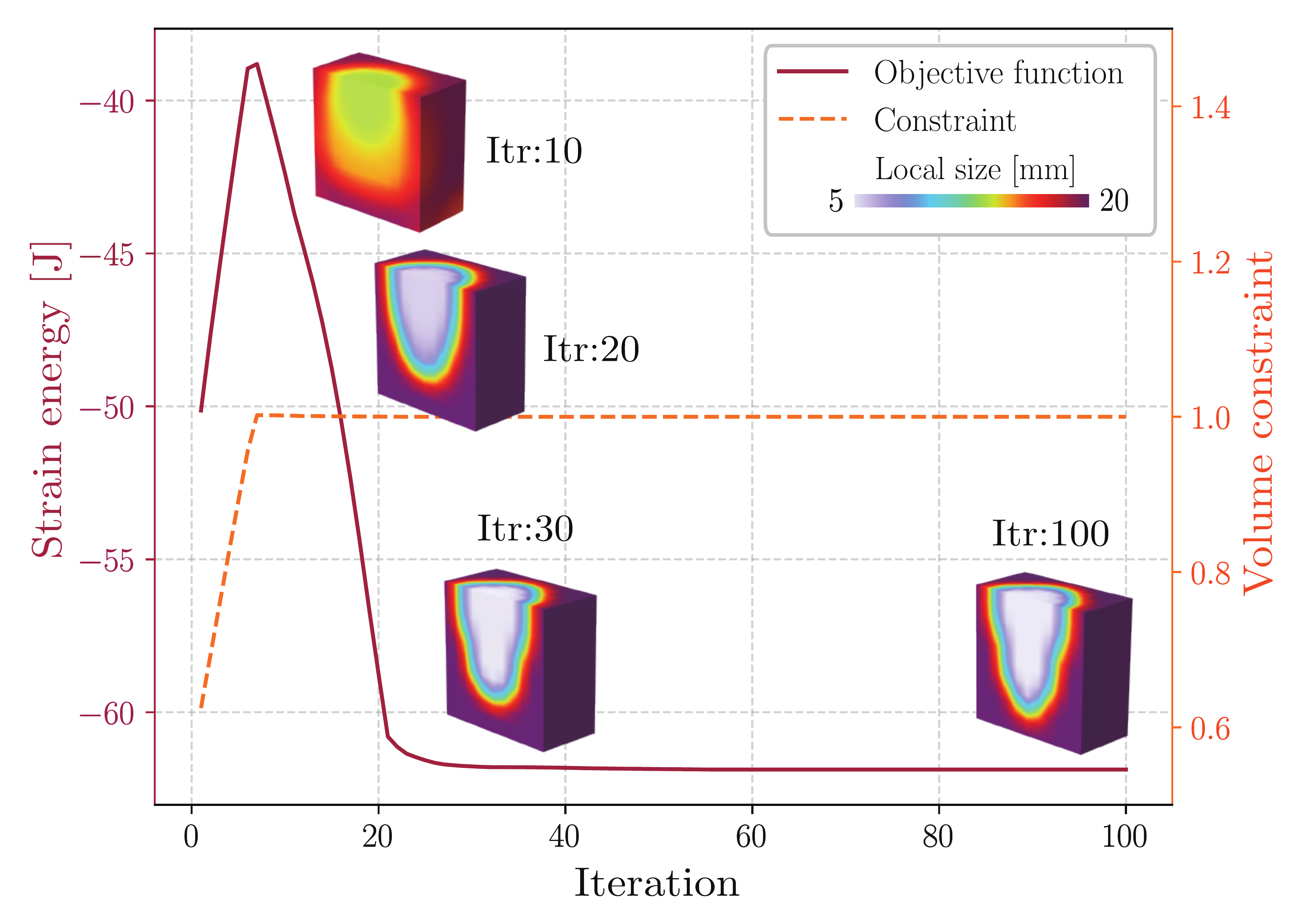}
  \caption{Optimization history: negative strain energy on the left y-axis and volume constraint $V/V_\text{max}$ on the right y-axis; the optimized size distribution cut at the center of the design domain at iteration 10, 20, 30, and 100.}
  \label{fig:optimization_history}
\end{figure}

\begin{figure*}[t]
\centering
\includegraphics[width=0.95\textwidth]{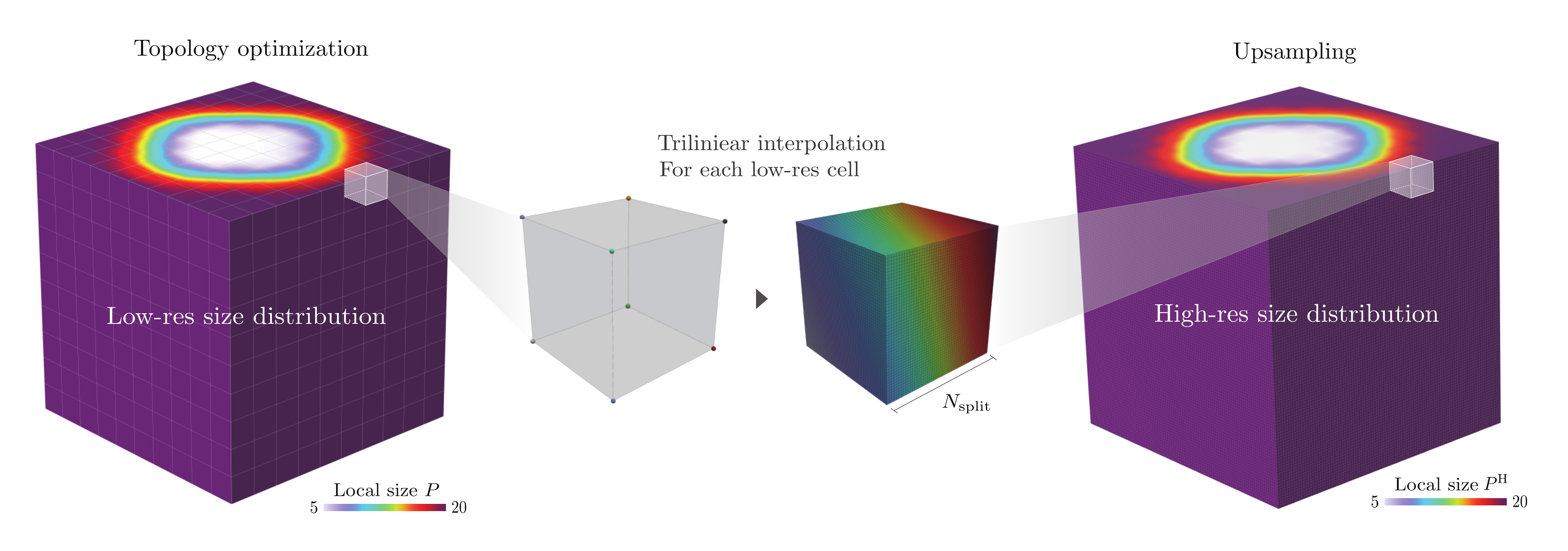}
\caption{Upsampling results and illustration of the upsampling process: the size optimization results from TO on the coarse mesh (left); the upsampled size optimization on the fine mesh (right); an example of the interpolation for a specific low-res cell (center).}
\label{fig:upsampling_process}
\end{figure*}

\begin{figure*}[t]
\centering
\includegraphics[width=0.95\textwidth]{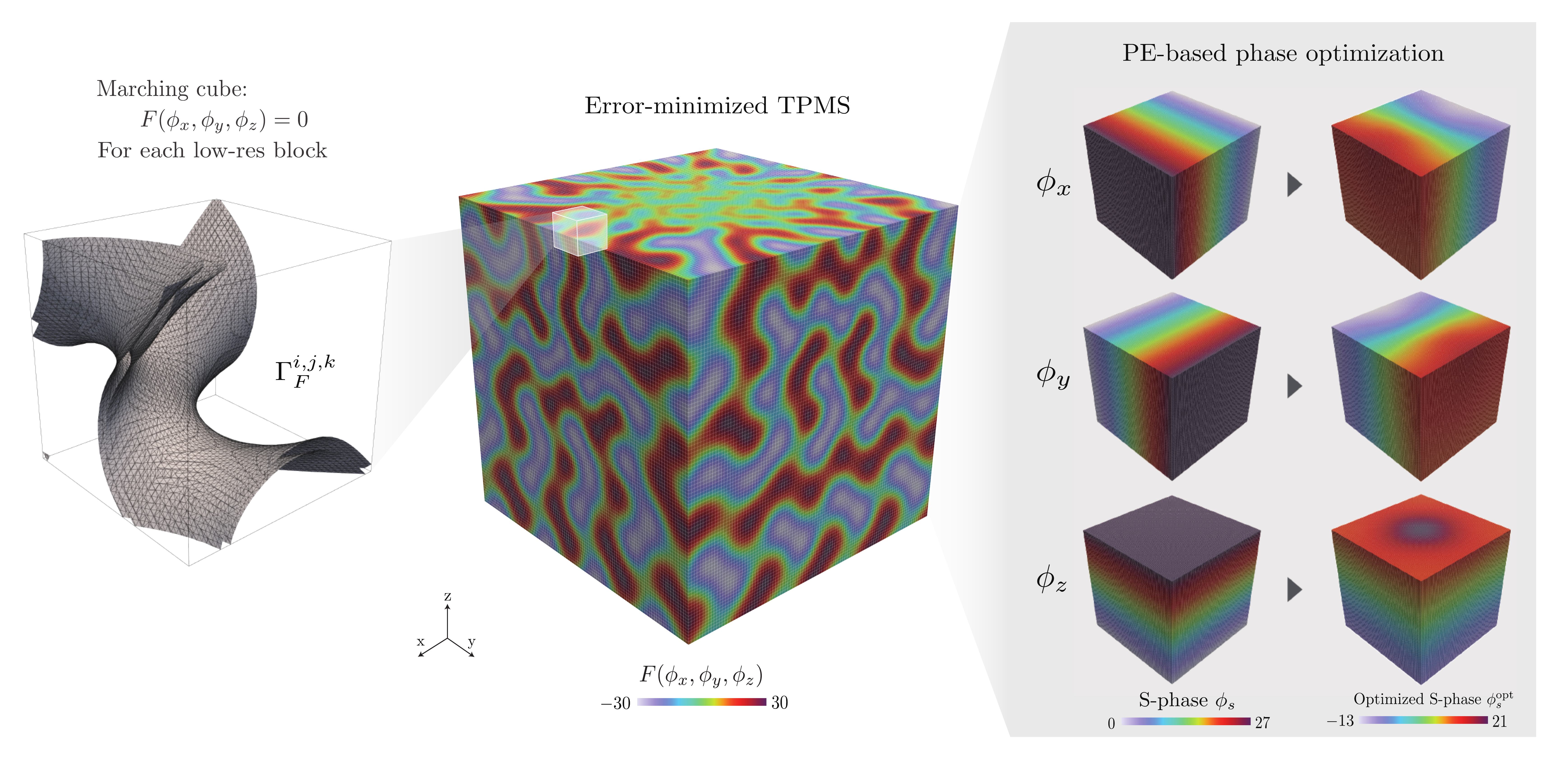}
\caption{Phase optimization results and marching cube illustration: the phase distributions $\phi_x$, $\phi_y$, and $\phi_z$ before and after the optimization (right); the TPMS function field $F(\phi_x, \phi_y, \phi_z)$ obtained from the optimized phase distributions (center); an example of the marching cube for a specific low-res cell (left).}
\label{fig:phase_opt_marching_cube}
\end{figure*}

\begin{figure*}[t]
\centering
\includegraphics[width=\textwidth]{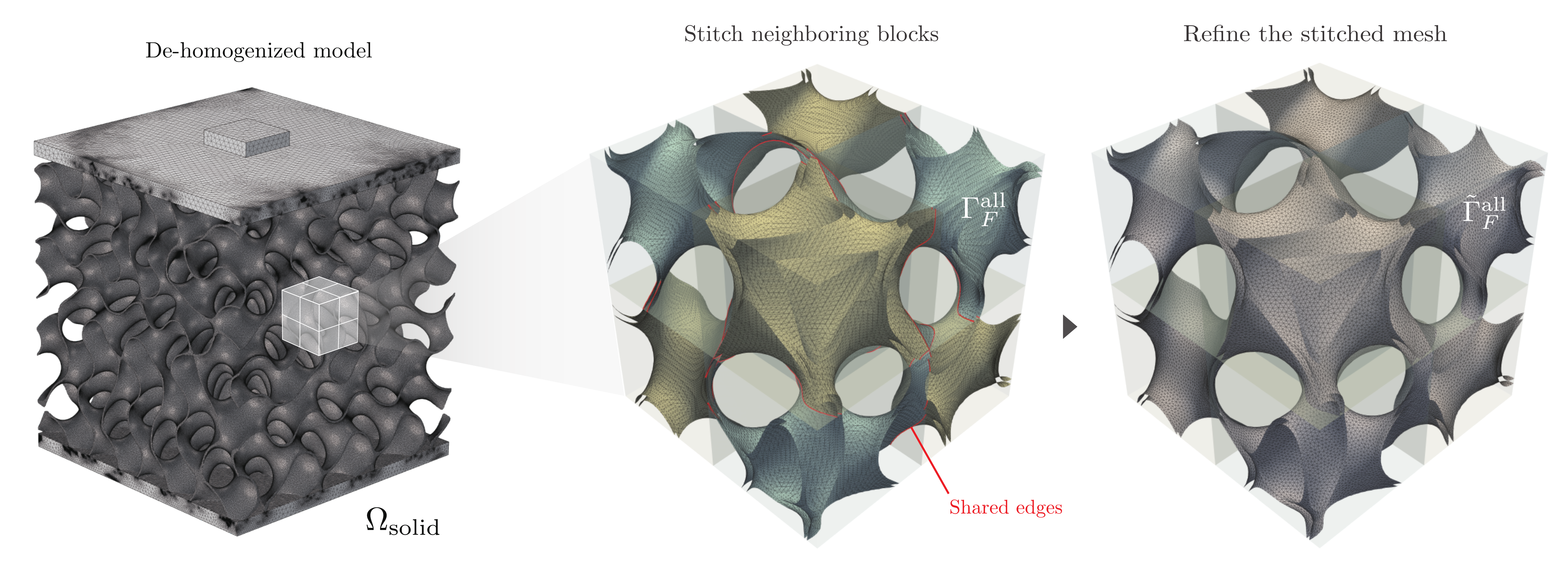}
\caption{Illustration of the stitching and remeshing process.}
\label{fig:stitch_remesh}
\end{figure*}

\begin{figure*}[t]
  \centering
  \includegraphics[width=0.9\textwidth]{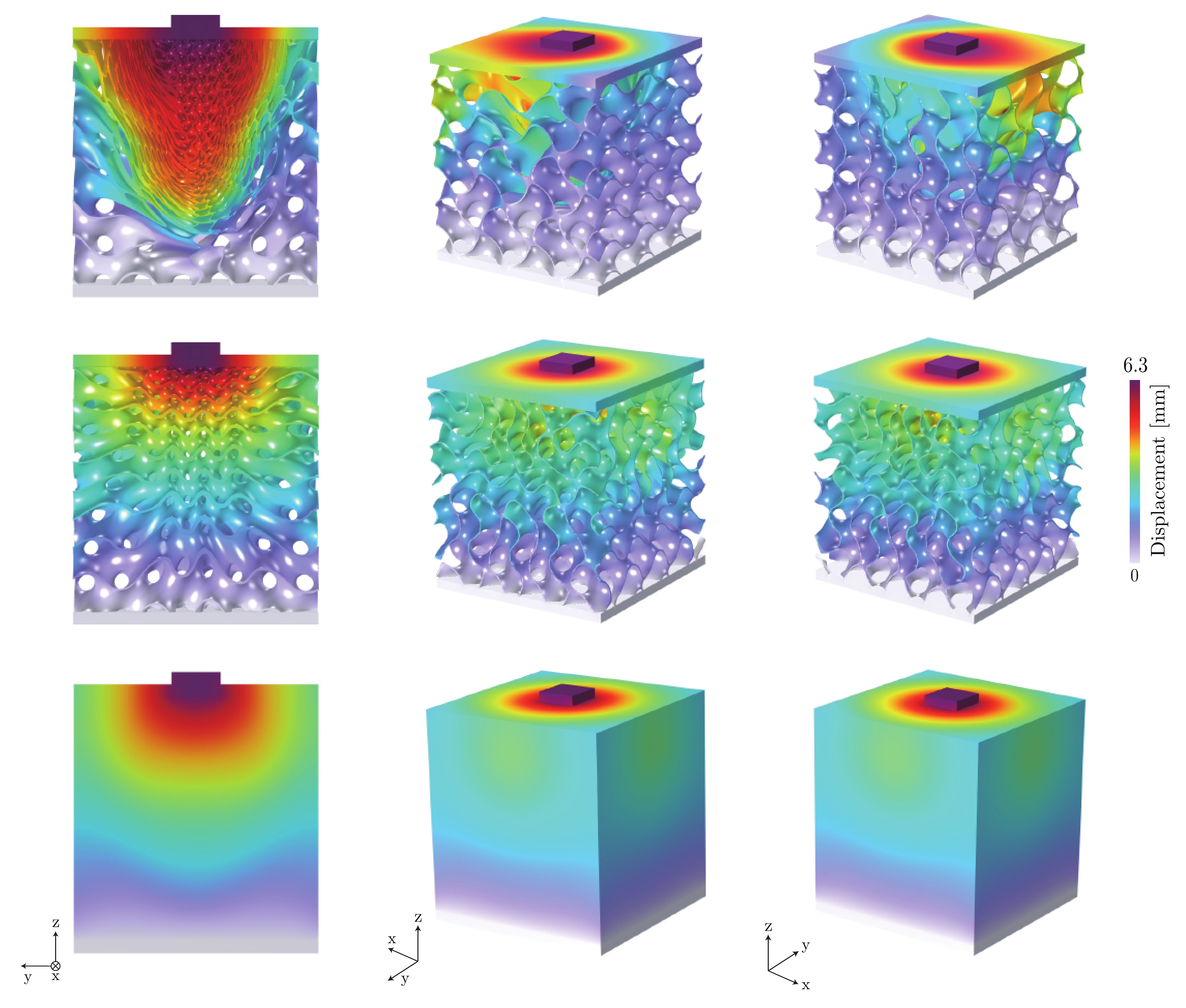}
  \caption{Displacement distribution of the models with the optimized size distribution obtained from the homogenized TO, showing the de-homogenized models from the PM method (top) and the proposed method (middle), as well as the homogenized model (bottom).}
  \label{fig:opt_disp}
\end{figure*}

\begin{figure*}[t]
  \centering
  \includegraphics[width=0.9\textwidth]{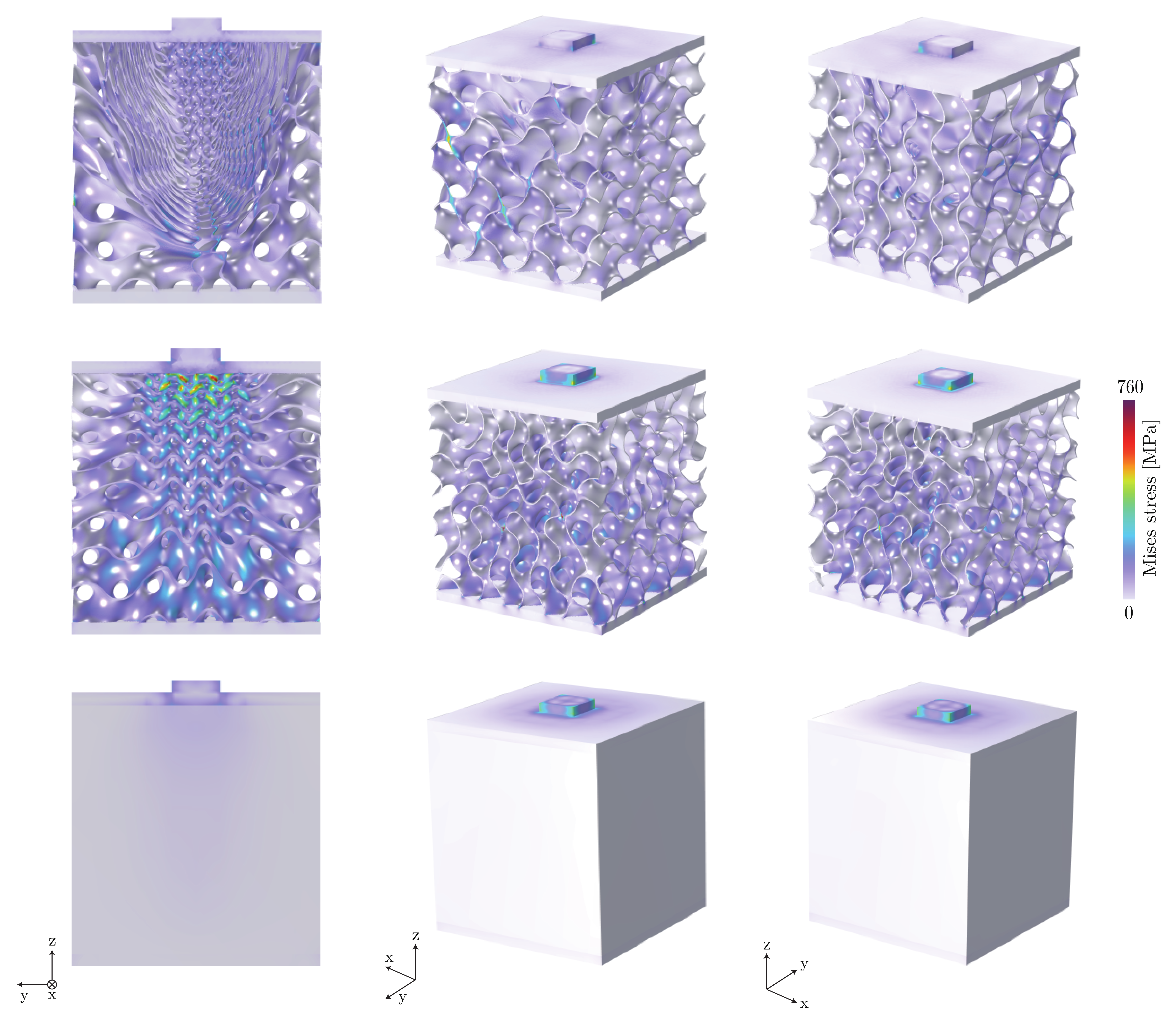}
  \caption{Stress distribution of the models with the optimized size distribution obtained from the homogenized TO, showing the de-homogenized models from the PM method (top) and the proposed method (middle), as well as the homogenized model (bottom).}
  \label{fig:opt_stress}
\end{figure*}

\begin{figure}[t]
  \centering
  \includegraphics[width=\columnwidth]{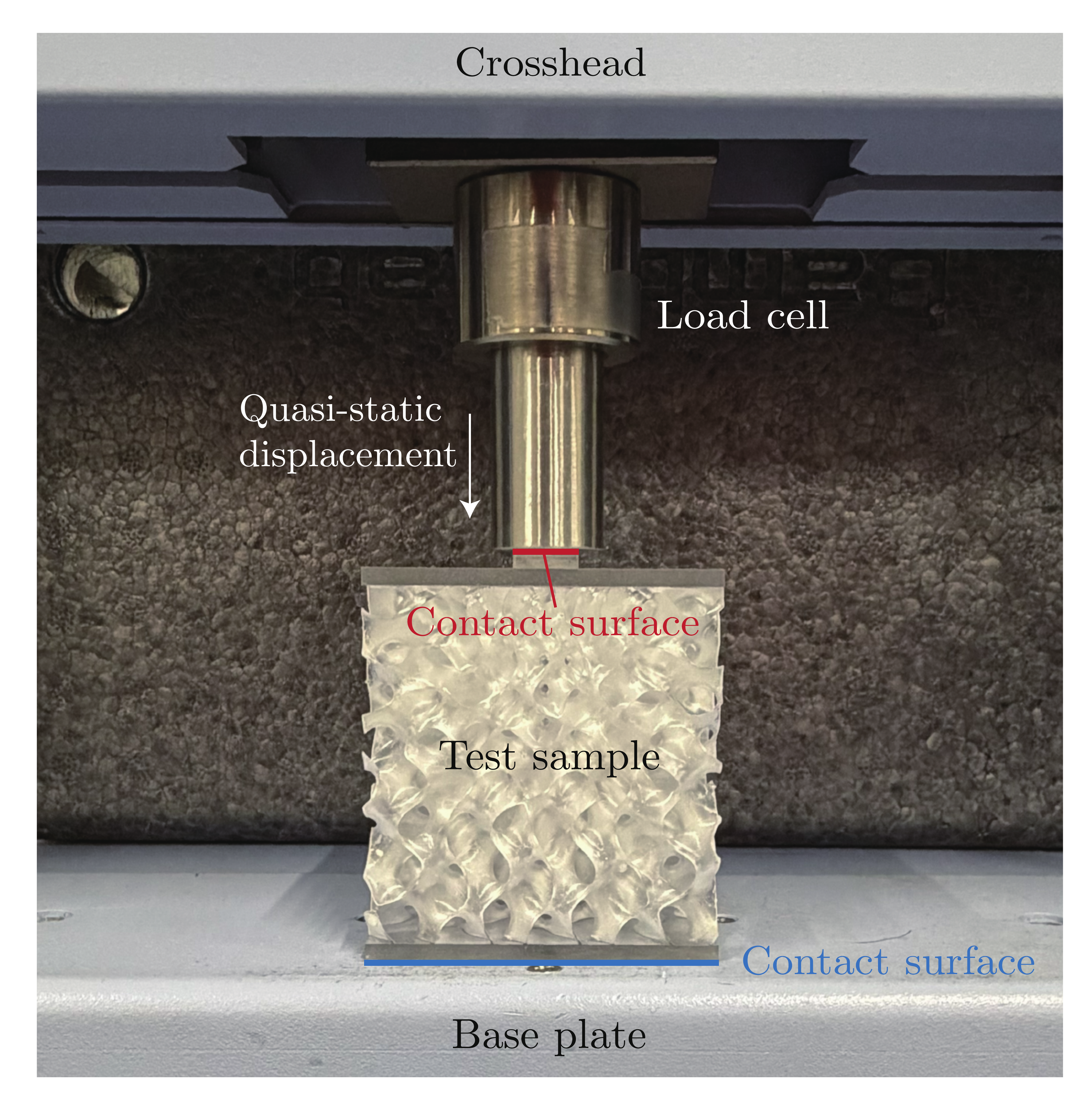}
  \caption{Experimental setup for the compression tests.}
  \label{fig:experiment_setup}
\end{figure}

\begin{figure}[t]
  \centering
  \includegraphics[width=\columnwidth]{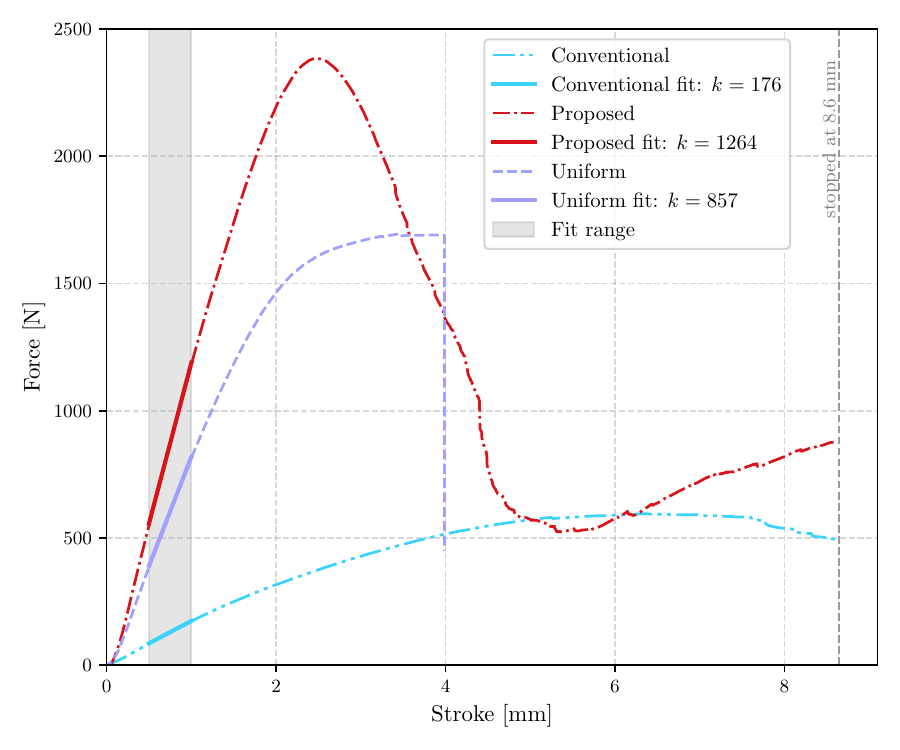}
  \caption{Load-displacement curves for the three results: uniform, optimized with the proposed method, and optimized with the conventional method. The solid regions with gray background represent the linear region that was used for effective stiffness calculation.}
  \label{fig:group_A_vs_B_vs_C_mean}
\end{figure}

\begin{figure*}[tb]
  \centering
  \includegraphics[width=\textwidth]{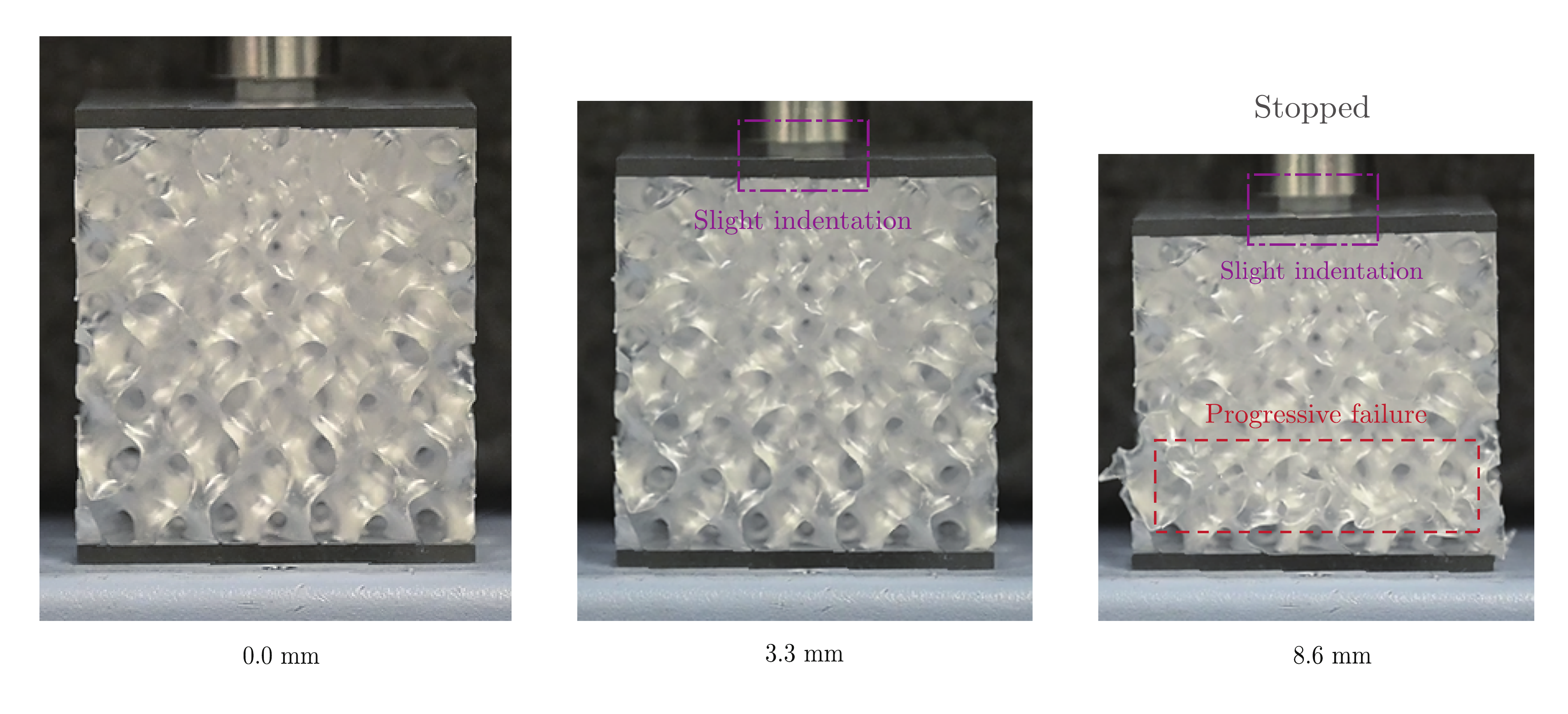}
  \caption{Side view of the specimens from the proposed method during the compression test. The bottom length value indicates the applied displacement for each image.}
  \label{fig:test_sideview_proposed}
\end{figure*}

\begin{figure*}[tb]
  \centering
  \includegraphics[width=\textwidth]{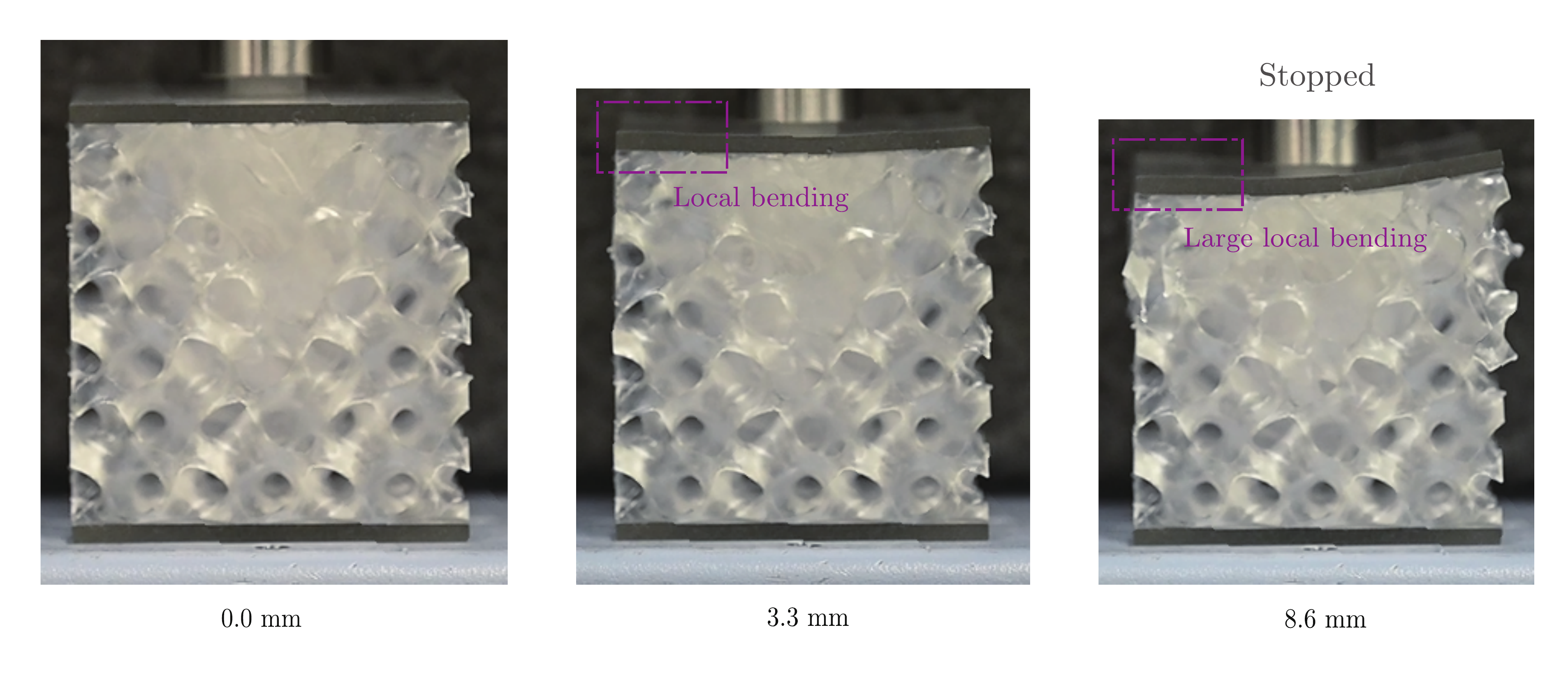}
  \caption{Side view of the specimens from the conventional PM method during the compression test. The bottom length value indicates the applied displacement for each image.}
  \label{fig:test_sideview_conventional}
\end{figure*}

\begin{figure*}[tb]
  \centering
  \includegraphics[width=\textwidth]{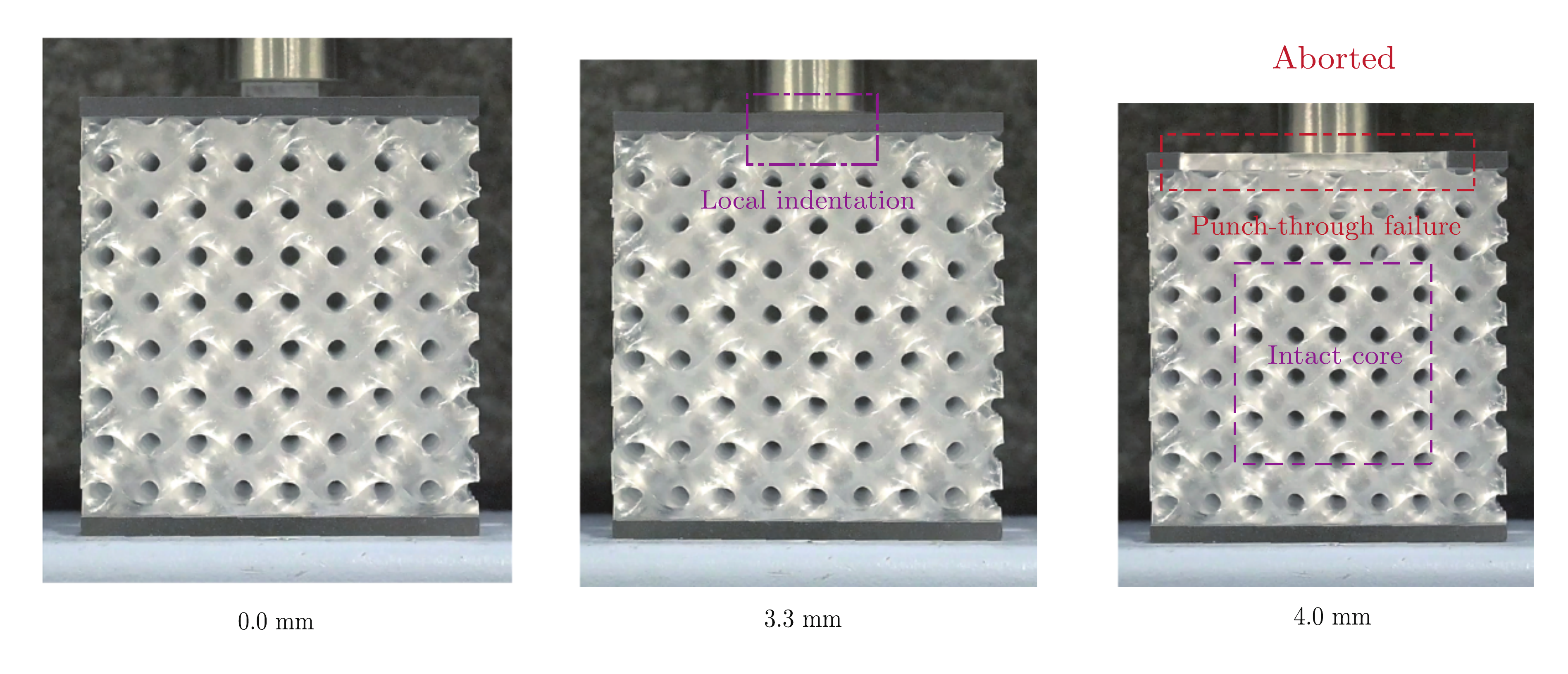}
  \caption{Side view of the specimens from the uniform case during the compression test. The bottom length value indicates the applied displacement for each image.}
  \label{fig:test_sideview_uniform}
\end{figure*}

\begin{figure*}[tb]
  \centering
  \includegraphics[width=0.75\textwidth]{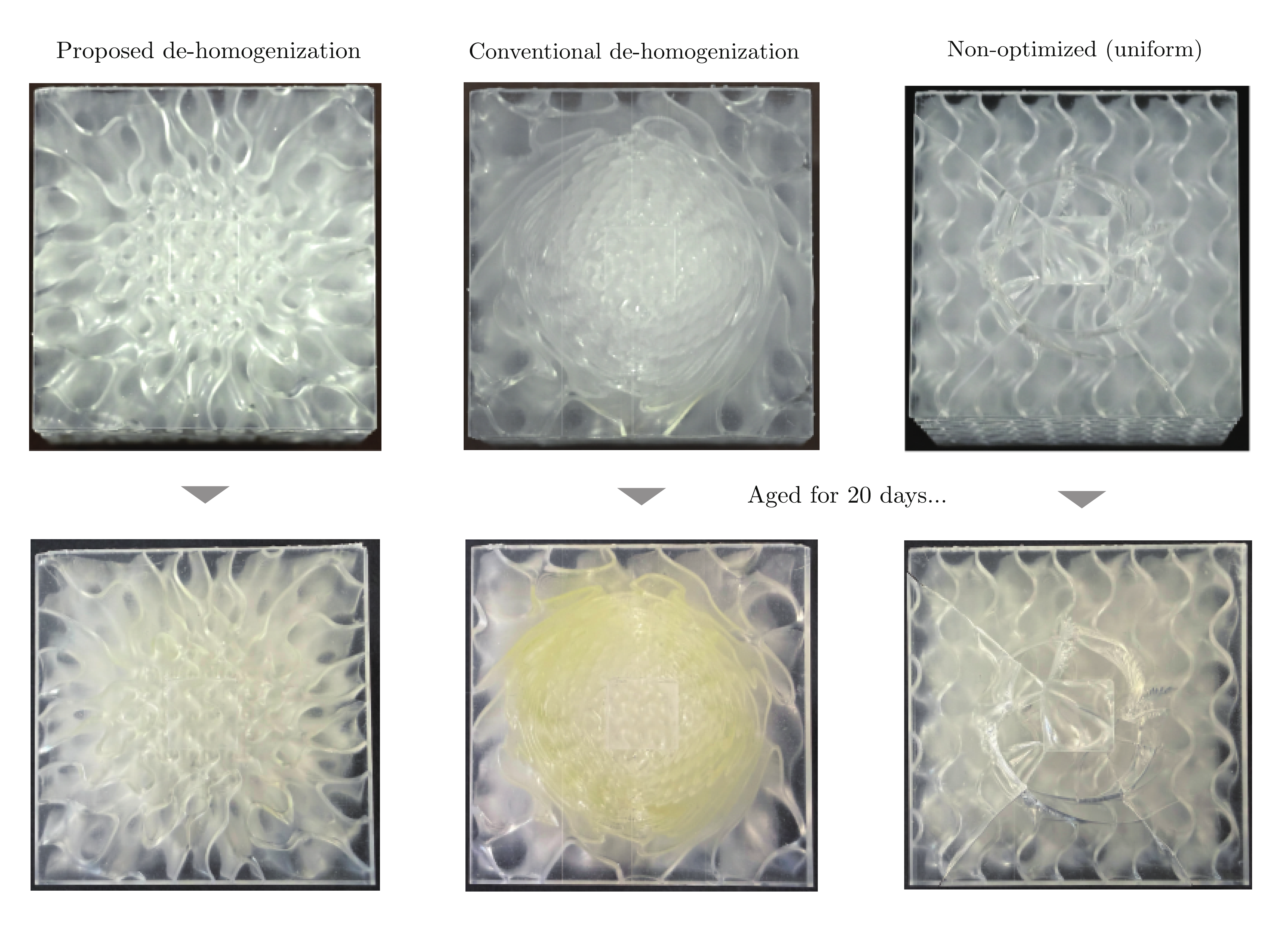}
  \caption{Top view of the test specimens before and 20 days after the compression test.}
  \label{fig:test_topview}
\end{figure*}

\subsection{Verification on homogenization}
\label{subsec:verification}
Before applying the proposed method to the optimized design, we first verify the homogenization process by comparing the structural responses of the homogenized and de-homogenized models with a uniform size distribution.
The uniform size was chosen to be the middle of the size range used in the optimization, $P=12.5$ [mm], to first confirm that the homogenization process accurately captures the structural behavior of the de-homogenized structure.
After the first verification, we further conducted the second verification for another uniform size, $P=14.8$ [mm], which corresponds to the same mass as the optimized design to be described in the next subsection.
This second verification facilitates the fair comparison between the uniform and optimized design with the same mass in terms of structural responses.
Table~\ref{tab:structural_response_comparison} summarizes the strain energy and mass of the homogenized and de-homogenized models for these two uniform size cases, where those metrics are calculated from the FEA under the same loading and boundary conditions as those used in the optimization.
Both cases show excellent agreement between the homogenized and de-homogenized models, with the errors in strain energy and mass being less than 5\% and 1.5\%, respectively.
This excellent agreement confirms the accuracy of the homogenization process in capturing the macroscopic structural behavior of the de-homogenized structure.

\begin{table}[bt]
\centering
\caption{Comparison of strain energy and mass between the homogenized and de-homogenized models for the uniform size cases.}
\label{tab:structural_response_comparison}
{\footnotesize
\renewcommand{\arraystretch}{1.15}
\begin{tabular}{llccc}
\hline
Size [mm] & Value & Strain energy [J] & Mass [g] \\
\hline
\multirow{3}{*}{12.5}
& Homo & 37.6 & 54.9 \\
& Dehomo & 36.3 & 54.5 \\
& Error & 3.6\% & 0.7\% \\
\hline
\multirow{3}{*}{14.8}
& Homo & 43.5 & 58.9 \\
& Dehomo & 41.6 & 58.2 \\
& Error & 4.7\% & 1.3\% \\
\hline
\end{tabular}
}
\end{table}

To further confirm the consistency between the homogenized and de-homogenized models, we also compare the displacement and stress distributions from the FEA for the second uniform size case, as shown in Fig.~\ref{fig:uniform_disp} and Fig.~\ref{fig:uniform_stress}, respectively.
As shown in Fig.~\ref{fig:uniform_disp}, in the de-homogenized model, the applied displacement on the top of the analysis domain radially propagates through the core lattice structures, as observed from the sectional and side views.
The homogenized model exhibits closely similar displacement patterns, indicating that the homogenized model accurately captures the overall deformation behavior of the de-homogenized structure.
However, the de-homogenized model shows local variations in displacement with asymmetry against the displacement direction, which is more prominent near the most outer regions of the lattice core, while the homogenized model shows a more symmetric displacement distribution.
This is because the de-homogenized model captures the geometric differences at the outer regions of the lattice core, where the local geometry highly affects the local stiffness under the free boundary conditions, while the homogenized model assumes a periodic unit cell and thus cannot capture such local geometric effects.
This limitation can be considered small enough for the homogenized model to be used for the optimization, as the overall displacement patterns closely resemble each other between the two models, and the error in the objective strain energy remains small as shown in Table~\ref{tab:structural_response_comparison}.

Unlike the displacement distribution, the stress distribution shows more significant differences between the homogenized and de-homogenized models, as shown in Fig.~\ref{fig:uniform_stress}.
The de-homogenized model exhibits more prominent stress concentrations, progressively distributed from the top to the bottom and from the inner to the outer regions of the lattice core.
This is physically reasonable as the applied displacement on the top of the analysis domain propagates through the lattice core and is eventually transmitted to the bottom, where the fixed boundary condition is applied.
The homogenized model, on the other hand, shows almost no stress concentration in the lattice core, except for the top plate, which is analyzed as a solid structure.
More specifically, the maximum von Mises stress in the homogenized model is 508 MPa, which is significantly smaller than 795 MPa in the de-homogenized model.
This is similar to the limitation observed in the displacement distribution, but it is more significant in the stress distribution as the local geometric effects have a more pronounced influence on the local stress concentrations.
This indicates that the homogenized model may not accurately capture the local stress concentrations in the de-homogenized structure, which makes it less suitable for stress-based optimization, compared to the stiffness-based optimization considered in this study.

\subsection{Optimized design}
\label{subsec:optimized}

Here, we present the results of the homogenization-based TO using the verified homogenized model, assuming the limitation in the stress distribution is acceptable for the stiffness-based optimization.
Fig.~\ref{fig:optimization_history} shows the optimization history of the proposed method, where the objective strain energy and volume constraints are plotted against the number of iterations, and the optimized size distributions are shown at representative iterations.
The strain energy is negative in this case as the optimization maximizes the stiffness of the structure, and thus the objective is set to minimize the negative strain energy.
The constraint value is calculated as the $V/V_\text{max}$, where $V$ is the current volume fraction initialized as $0.5$ and $V_\text{max}$ is the maximum allowed volume fraction, i.e., $0.8$.
The iteration history shows that the optimization converges around 50 iterations, where the objective strain energy is significantly increased (i.e., the stiffness is improved under the fixed displacement) while the volume constraint is satisfied after the initial iterations.
The optimized size distribution shows a clear trend of smaller sizes near the top and larger sizes near the bottom and outer regions of the lattice core.
This is physically reasonable as the smaller sizes provide higher stiffness to resist the applied displacement on the top, while the larger sizes help to reduce the mass of the structure.

After the TO on the low-resolution homogenized model, the de-homogenization process is performed to generate the full-scale de-homogenized TPMS structure, by the following steps.
\begin{itemize}
\item \textbf{Upsampling}: The optimized size distribution from the low-resolution homogenized mesh is upsampled to a high-resolution size distribution $P^\text{H}(\mb{r})$ on the high-resolution mesh, using a trilinear interpolation~\cite{bolch2003marching} for each low-resolution cell.
Each low-resolution cell is split into $N_\text{split}^3$ high-resolution cells, where $N_\text{split}=130$ is used in this study.
Fig.~\ref{fig:upsampling_process} shows the upsampled size distribution on the high-resolution mesh, where the size distribution is smoothly interpolated from the low-resolution mesh, so that the overall trend of the optimized size distribution is preserved in the upsampled size distribution.
Note that the size distribution $P(\mb{r})$ in the de-homogenization formulation above is replaced with the upsampled size distribution $P^\text{H}(\mb{r})$ for the subsequent phase optimization and marching cube steps, in this homogenization-based TO framework.
\item \textbf{Phase optimization}: The PE-based phase optimization is performed to minimize the discrepancy between the target wavenumber fields obtained from the upsampled size distribution $P^\text{H}(\mb{r})$ and the actual wavenumber fields.
Fig.~\ref{fig:phase_opt_marching_cube} shows the optimized phase distribution and the corresponding TPMS function $F(\phi_x, \phi_y, \phi_z)$ calculated from the optimized phase distributions.
The optimized phase distribution exhibits significant deviations from the initial phase distribution, which is used for the conventional PM method, implying that the phase optimization is necessary to achieve the optimized structural performance after de-homogenization.
\item \textbf{Marching cubes}: The marching cubes algorithm is applied to each low-resolution block that contains the high-resolution mesh with the optimized size and phase distributions, to extract the iso-surface of the TPMS structure $\Gamma_F$.
The intention of applying the marching cubes to each block is to avoid the memory issue of GPUs, as the marching cubes for the entire high-resolution mesh often requires more memory than the available GPU memory.
\item \textbf{Stitching}: The extracted TPMS structures from each block are stitched together at the shared edges to form the full-scale de-homogenized mesh $\Gamma_F^\text{all}$, as shown in Fig.~\ref{fig:stitch_remesh}.
The continuity of the structure across the block boundaries is ensured by first optimizing the phase distribution in the entire design domain and then applying the marching cubes algorithm to each block, which allows for the exact matching of the edge geometry across the block boundaries.
\item \textbf{Remeshing}: The meshes $\Gamma_F^\text{all}$ obtained from the marching cubes are often of low quality, with irregularly shaped elements and non-uniform element sizes, which can lead to inaccurate FEA results and convergence issues, as shown in Fig.~\ref{fig:stitch_remesh}. Therefore, the stitched iso-surface is remeshed to generate a high-quality mesh $\tilde{\Gamma}_F^\text{all}$ for the FEA, which is necessary for accurately evaluating the structural performance of the de-homogenized structure.
\end{itemize}

\begin{table}[bt]
\centering
\caption{Comparison of strain energy and mass between the homogenized and de-homogenized models obtained from the PM method and the proposed method.}
\label{tab:optimized_structural_response_comparison}
{\footnotesize
\renewcommand{\arraystretch}{1.15}
\begin{tabular}{llccc}
\hline
Method & Value/Error & Strain energy [J] & Mass [g] \\
\hline
Homo & Value & 61.9 & 58.9 \\
\hline
\multirow{3}{*}{Dehomo (PM)}
& Value & 22.5 & 57.8 \\
& Error & 63.6\% & 1.8\% \\
& vs Uniform & -45.8\% & -0.5\% \\
\hline
\multirow{3}{*}{Dehomo (Proposed)}
& Value & 62.4 & 58.4 \\
& Error & 0.8\% & 0.9\% \\
& vs Uniform & +50.1\% & +0.4\% \\
\hline
\end{tabular}
}
\end{table}

To demonstrate the effectiveness of the proposed PE-based phase optimization, we also compare the de-homogenized structure obtained from the conventional PM method using Gaussian smoothing with $\sigma_c=N_\text{split}$, with that from the proposed method.
Note that the same de-homogenization process except for the phase optimization is applied to the uniform case shown in the verification and PM case.

Table~\ref{tab:optimized_structural_response_comparison} summarizes the strain energy and mass of the homogenized and de-homogenized models obtained from the optimized design, including the error of the de-homogenized models compared to the homogenized model, and the comparison with the uniform case with almost the same mass.
Both the PM method and the proposed method maintain the mass of the de-homogenized structure close to that of the homogenized model, with the errors being less than 2\%.
However, the PM method results in a significant decrease in the strain energy compared to the homogenized model, with an error of 63.6\%, while the proposed method achieves a much smaller error of 0.8\%.
As a result, the proposed method achieves a significant improvement in the strain energy compared to the uniform case with an increase of 50.1\%, as optimized in the homogenized model, while the PM method shows a decrease of 45.8\%, completely failing to achieve the intended improvement after de-homogenization.
This highlights the importance of the phase optimization in the de-homogenization process, as it allows for the generation of a de-homogenized structure that can achieve the desired structural performance after de-homogenization, which is otherwise not guaranteed by the conventional PM method.

Fig.~\ref{fig:opt_disp} and Fig.~\ref{fig:opt_stress} show the displacement and stress distributions of the models with the optimized size distribution, respectively, showing the de-homogenized models from the PM method and the proposed method, as well as the homogenized model.
As shown in both figures, the PM method with Gaussian smoothing results in extremely distorted structures particularly around the regions with intense size changes, whereas the proposed method successfully generates well-preserved Gyroid structures that much more closely reflect the optimized size distribution.
Consequently, the PM method shows significantly different displacement distribution from the homogenized model, with the applied displacement transmitting through the distorted lattice core beyond the central region.
At the same time, the size distribution in the PM method exhibits a separate pattern between the inside and outside of the distorted region, with inside the distorted region deforming almost independently from the outside with small displacement gradients.
Since the local strain energy is calculated as the product of the local stress and strain, the local strain energy in the distorted region is significantly reduced due to the small displacement gradients, which leads to the significant decrease in the overall strain energy of the structure, as shown in Table~\ref{tab:optimized_structural_response_comparison}.
Additionally, the PM method shows significant asymmetry in the displacement distribution against the displacement direction, which is not observed in the homogenized model.
This is because the local distortion of the lattice core is highly dependent on the distance from the origin, which is here a corner of the analysis domain, i.e., the bottom-right corner in the sectional view.

In contrast, the proposed method shows closely similar displacement patterns to the homogenized model.
The applied displacement on the top more entirely propagates through the optimized lattice core, compared to the uniform case in the verification.
This more progressive propagation of the applied displacement can explain the significant improvement in the strain energy compared to the uniform case, as the higher displacement gradients in the larger regions can lead to higher local strain energy.

Similar to the uniform case, the stress distribution shows more significant differences between the homogenized and de-homogenized models, compared to the displacement distribution, for both the PM method and the proposed method.
However, the optimized homogenized model shows larger stress distribution in the lattice core compared to the uniform case, which can also explain the significant improvement in the strain energy compared to the uniform case, as the local strain energy is calculated as the product of the local stress and strain.
The de-homogenized model from the proposed method exhibits a larger stress distribution in more global regions of the lattice core, compared to the uniform case, which also aligns with the improvement in the strain energy.
The de-homogenized model from the PM method, on the other hand, shows smaller stress distribution entirely, which is consistent with the significantly smaller strain energy than that of the homogenized model, as shown in Table~\ref{tab:optimized_structural_response_comparison}.

\section{Experimental validation}
\label{sec:exp}

To validate the effectiveness of the proposed method, we conducted quasi-static compression tests on the de-homogenized models obtained from the proposed method and the conventional PM method, as well as the uniform case.
After the tests, we compared their structural responses in terms of load-displacement curves and failure modes.

\subsection{Experimental setup}
\label{subsec:setup}

The compression tests were conducted using a universal testing machine, SHIMADZU AGS-X, equipped with a calibrated load cell with a maximum capacity of 10 kN, as shown in Fig.~\ref{fig:experiment_setup}.
The test specimens were placed between the base plate on the bottom and the loading cylinder on the top, with the loading direction aligned with the vertical axis.
No lateral support was provided, allowing free deformation in the horizontal direction, which technically differs from the fixed boundary conditions assumed in the FEA.
This setup, however, can be considered as a more realistic representation of the loading conditions, since the applied load was vertical and the lateral deformation was observed to be small enough to not affect the overall structural response.

The tests were performed on three samples for each of the three cases: uniform, optimized with the proposed method, and optimized with the conventional PM method, resulting in a total of nine test specimens.
The test specimens were fabricated using a stereolithography (SLA) 3D printer, Formlabs Form 4, with a layer thickness of 0.1 mm.
The material used for printing was the same resin material as used in the FEA.
The fabrication process consisted of the following steps: printing the specimens aligned with the $y$ direction from negative to positive, washing in tripropylene glycol monomethyl ether using Form Wash for 10 minutes, and curing under UV light for five minutes at room temperature using Form Cure, according to the instructions provided by the printer manufacturer.

For one specimen of each case, the specimens were compressed at a constant displacement rate of 0.5 mm/min until failure or a maximum displacement of 8.6 mm, which corresponds to 12\% strain, were reached, whichever occurred first.
After confirming the failure mode and the load-displacement curve including the non-linear region for one specimen of each case, the test was conducted for the remaining two specimens of each case at the same displacement rate, with the test stopped a while after the end of the linear region.
The stopping points were determined for each case based on the load-displacement curve of the first specimens.
This is because the homogenized and de-homogenized models assume linear elasticity, making it more fair to compare the structural responses mainly in the linear region.
The non-linear region can be affected by various factors such as local buckling and material non-linearity, which are not captured in the homogenized and de-homogenized models.

\subsection{Experimental results}
\label{subsec:results}

To investigate the fabrication accuracy of the de-homogenized models, we measured the mass of the fabricated specimens and compared them with the mass calculated from the de-homogenized models, as shown in Table~\ref{tab:mass_comparison}.
The error is calculated from the mean mass of the three specimens for each case.
The results show that the fabricated specimens from the proposed method and the uniform case have a similar mass with a small error of around 10\%, while the specimens from the PM method have a significantly larger mass with a large error of around 30\%.
This large mass error is due to the resin trapped in the distorted regions, which is further confirmed in the observation of the specimens after the tests as discussed at the end of this section.
This indicates that the PM method makes it difficult to accurately fabricate the optimized structures, due to the strong distortion of the lattice structures.

\begin{table}[t]
  \centering
  \caption{Comparison of the mass between the fabricated specimens and the de-homogenized models.}
  \label{tab:mass_comparison}
  {
    \footnotesize
    \renewcommand{\arraystretch}{1.1}
    \begin{tabular}{lccc}
    \hline
    Case & Specimen [g] & Dehomo model [g] & Error [\%] \\
    \hline
    Uniform & 63.4 $\pm 0.22$ & 58.2 & 9.0 \\
    Proposed method & 64.4 $\pm$ 0.87 & 58.4 & 10 \\
    PM method & 74.0 $\pm$ 0.87 & 57.8 & 28 \\
    \hline
    \end{tabular}
  }  
\end{table}

Fig.~\ref{fig:group_A_vs_B_vs_C_mean} shows the load-displacement curves for the first specimen of each case.
The slopes of the load-displacement curve are calculated in the linear region for all three specimens for each case, to measure the effective stiffness of the specimens, as summarized in Table~\ref{tab:effective_stiffness_comparison}.
The table and the curves in Fig.~\ref{fig:group_A_vs_B_vs_C_mean} show that the proposed method remarkably improves the effective stiffness by averagely 54.2\% compared to the uniform case, which is consistent with the 50.1\% improvement in the FEA.
In contrast, the PM method shows a significant decrease in the effective stiffness by averagely 77.3\% compared to the uniform case, which is more significant than the 45.8\% decrease in the FEA.
This significant decrease in the effective stiffness can be attributed to the strong distortion of the lattice structures and the consequent fabrication inaccuracy in the PM method.

\begin{table}[t]
  \centering
  \caption{Comparison of the effective stiffness between the fabricated specimens and the homogenized models.}
  \label{tab:effective_stiffness_comparison}
  {
    \footnotesize
    \renewcommand{\arraystretch}{1.1}
    \begin{tabular}{lcc}
    \hline
    Case & Slope [N/mm] & vs Uniform [\%]\\
    \hline
    Uniform & 827.6 $\pm$ 27.0 & - \\
    Proposed method & 1276.1 $\pm$ 12.4 & +54.2 \\
    PM method & 188.1 $\pm$ 10.1 & -77.3 \\
    \hline
    \end{tabular}
  }
\end{table}

In the load-displacement curves, only the specimens from the uniform case show a clear failure with a sudden drop in the load, whereas the specimens from the proposed method and the PM method kept sustaining the load without any clear failures until we stopped the test.
To analyze the failure modes, we observed the specimens during the test and after the test, where the side view from the proposed method, the PM method, and the uniform case are shown in Fig.~\ref{fig:test_sideview_proposed}, Fig.~\ref{fig:test_sideview_conventional}, and Fig.~\ref{fig:test_sideview_uniform}, respectively, while the top view of all cases is shown in Fig.~\ref{fig:test_topview}.
The specimen from the uniform case shows a local indentation on the loading surface at the displacement of 3.3 mm, and then the button-like box structure on the top suddenly collapses punching through the top surface at the displacement of 4.0 mm, with the inner lattice core remaining intact.
This failure mode can be observed from the top view in Fig.~\ref{fig:test_topview}, where the uniform case specimen shows a clear circular crack and indentation around the loading surface.
The specimens from the proposed method and the PM method, on the other hand, show no large failures during the test, but gradually deformed with small local failures, as observed from the side view in Fig.~\ref{fig:test_sideview_proposed} and Fig.~\ref{fig:test_sideview_conventional}, respectively.
The proposed method specimen shows slight indentation on the loading surface at the displacement of 3.3 mm, which is consistent with the displacement distribution in the FEA, and at the end of the test, the specimen shows progressive deformation and failure with densification of the lattice core.
This densification can explain the secondary increase in the load after around the displacement of 5 mm, as the densification can lead to the increase in the local stiffness of the lattice core.
Contrarily, the PM method specimen shows local bending of the top plate at the left corner at the displacement of 3.3 mm, which became more pronounced at the end of the test, with obviously asymmetric deformation with respect to the loading direction.
This deformation closely matches the asymmetric displacement distribution in the FEA, as shown in Fig.~\ref{fig:opt_disp}.

Moreover, we observed that after 20 days from the tests, the specimens from the proposed method showed very slightly yellowed color around the most dense regions of the lattice core, and the specimens from the PM method showed even more obviously yellowed color around the same regions.
This implies that the yellowing can be caused by the remaining resin trapped inside the dense regions, confirming the mass comparison results in Table~\ref{tab:mass_comparison}, where the PM method specimens have a significantly larger mass than the designed mass, and the proposed method specimens have a slightly larger mass.

\section{Conclusions}
\label{sec:con}

In this study, we proposed a homogenized topology optimization (TO) method for the spatial cell-size optimization of triply-periodic minimal surface (TPMS) lattice structures, which can accurately achieve the optimized structural response after de-homogenization.
To achieve this, we introduced a novel de-homogenization technique that directly minimizes the difference between the wavenumber fields obtained from the target and actual size distributions.
This problem can be efficiently solved as a typical Poisson's equation using the discrete cosine transform (DCT).
We first verified the proposed de-homogenization method through numerical examples, showcasing its capability in significantly reducing the known distortion of the de-homogenized TPMS structures from the conventional periodic modulation (PM) method.
After that, we applied the proposed method to a numerical example of a stiffness maximization problem under fixed displacement.
In this example, we compared the optimized structural response of the de-homogenized models from the proposed method and the PM method, with that of the homogenized model and a reference model with uniform cell size.
Finally, we fabricated the optimized structures and the reference model using a stereolithography 3D printer, and conducted quasi-static compression tests to validate the effectiveness of the proposed method in improving the structural response compared to the conventional PM method.

Here are the key findings from the numerical examples and experimental validation:
\begin{itemize}
    \item The numerical results show that the proposed method successfully reduced the distortion of the de-homogenized structures, leading to an excellent agreement by only 0.8\% in the strain energy between the homogenized and de-homogenized models, as opposed to 63.6\% difference in the conventional PM method.
    \item The reduced distortion from the proposed method allows the optimized structure to significantly improve the strain energy by 50.1\% compared to the uniform case with the same mass on the de-homogenized models, while the PM method results in a significant decrease by 45.8\%.
    \item The experimental validation shows that the effective stiffness of the optimized structure from the proposed method is 54.2\% higher than that of the uniform case, while the PM method results in a significant decrease by 77.3\% even though exactly the same size distribution was used for the de-homogenization in both methods.
\end{itemize}

Despite its effectiveness in improving the de-homogenized structural response, the proposed method still has some limitations.
First, the homogenized model is currently unable to capture the local asymmetry in the displacement distribution observed in the de-homogenized model, due to the isotropic and periodic assumption in the homogenized model.
Second, the stress distribution shows significant differences between the homogenized and de-homogenized models, for both the PM method and the proposed method, due to the macroscopic nature of the homogenized model.
Future work could focus on developing a more accurate homogenized model that can overcome these limitations, as well as exploring other optimization objectives such as strength and durability, which are more directly related to the stress distribution.

\section*{Declaration of competing interest}
The authors declare that they have no known competing financial interests or personal relationships that could have appeared to influence the work reported in this paper.

\section*{Acknowledgements}
This work was partly supported by Innovative Future Space Transportation Systems Research and Development Program, Japan Aerospace Exploration Agency (JAXA). 
This work was also supported by JSPS KAKENHI Grant Number 23H01323.
The authors gratefully acknowledge Dr. Sangmin Lee for valuable guidance and assistance in conducting the experiments at F3D Lab, Sanken, The University of Osaka.

\appendix

\section{Thickness specification for the solid wall}
\label{ap:thickness}

\begin{figure}[ht]
  \centering
  \includegraphics[width=\columnwidth]{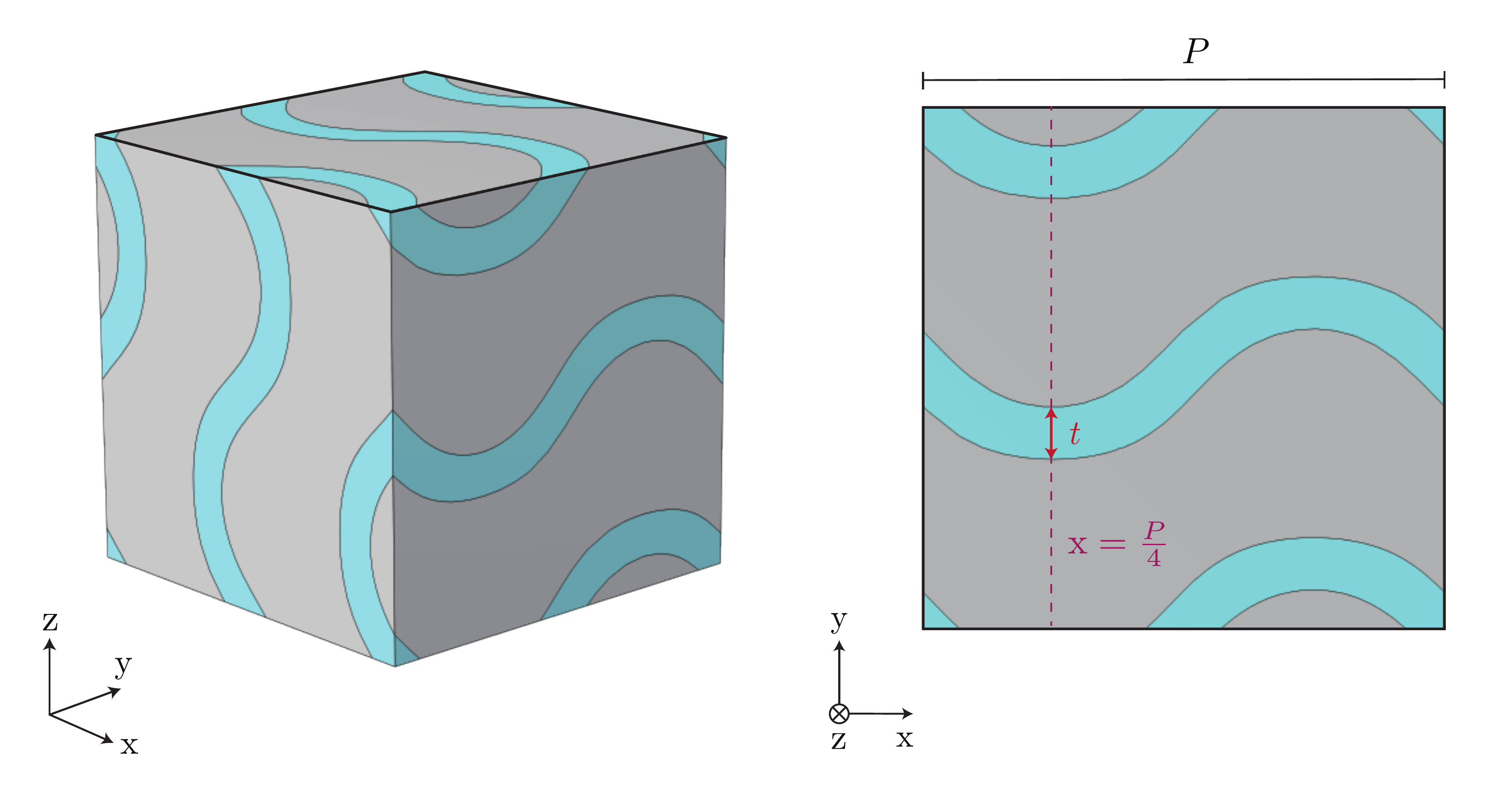}
  \caption{Thickness specification for the solid wall. Light blue regions represent the solid wall, and the purple dashed line represents the line of interest.}
  \label{fig:thickness_specification}
\end{figure}

When the $c$ value is fixed throughout the unit cell, the thickness of the solid wall is not uniform across the unit cell.
This is because when using Eq.~\eqref{eq:tpms_solid} to generate the solid wall, the local thickness is essentially determined by the local distance between two offset iso-surfaces of the implicit function, which is not constant across the unit cell.
Therefore, to achieve a thickness specification for the solid wall at reasonable cost, we can determine the $c$ value based on the desired thickness at a specific location, such as the location at quarter of the unit cell length for Gyroid, as shown in Fig.~\ref{fig:thickness_specification}.
Substituting the phase function from the PM method in Eq.~\eqref{eq:phase_traditional} into Eq.~\eqref{eq:type}, the TPMS function can be expressed as follows:
\begin{equation}
  \begin{aligned}
    F_G (\mb{r}) & = \sin \left(\frac{2\pi}{P(\mb{r})} x \right) \cos \left(\frac{2\pi}{P(\mb{r})} y \right) + \sin \left(\frac{2\pi}{P(\mb{r})} y \right) \cos \left(\frac{2\pi}{P(\mb{r})} z \right) \\
      & \ + \sin \left(\frac{2\pi}{P(\mb{r})} z \right) \cos \left(\frac{2\pi}{P(\mb{r})} x \right) = 0.  
  \end{aligned}
\end{equation}
\noindent Then, the TPMS function can be rearranged on the line of $x=P/4$ on the top surface of the unit cell, as follows:
\begin{equation}
  \begin{aligned}
    F_G(P/4, y, P) & = \sqrt{2}\sin \left(\frac{2\pi}{P} y + \frac{3\pi}{4} \right).
  \end{aligned}
\end{equation}
\noindent Substituting this into Eq.~\eqref{eq:tpms_solid}, the range of the $y$ can be expressed as follows:
\begin{equation}
  \frac{P}{2\pi} \left\{ \sin^{-1} \left( -\frac{c}{\sqrt{2} P} \right) - \frac{3\pi}{4} \right\} \leq y \leq \frac{P}{2\pi} \left\{ \sin^{-1} \left( \frac{c}{\sqrt{2} P} \right) - \frac{3\pi}{4} \right\}, \\
\end{equation}
\noindent and thus the thickness of the solid wall at this location can be expressed as follows:
\begin{equation}
  t = \frac{P}{2\pi} \left\{ \sin^{-1} \left( \frac{c}{\sqrt{2} P} \right) - \sin^{-1} \left( -\frac{c}{\sqrt{2} P} \right) \right\} = \frac{P}{\pi} \sin^{-1} \left( \frac{c}{\sqrt{2} P} \right).
\end{equation}
\noindent Rearranging this equation, the $c$ value can be determined based on the desired thickness $t$ as shown in Eq.~\eqref{eq:thickness_specification}.

  \bibliographystyle{elsarticle-num} 
  \bibliography{reference}





\end{document}